\documentclass{article}
\usepackage{amsfonts}
\usepackage{amsmath}
\usepackage{amssymb}
\usepackage{latexsym}

\newtheorem{theorem}{Theorem}

\newtheorem{corollary}[theorem]{Corollary}

\newtheorem{proposition}[theorem]{Proposition}
\newtheorem{remark}[theorem]{Remark}

\begin{document}

\author{George Avalos\thanks{%
email address: gavalos@math.unl.edu. The research of G. Avalos was partially
supported by the NSF grants DMS-1211232 and DMS-1616425.} \\
Department of Mathematics\\
University of Nebraska-Lincoln, 68588 USA \and Pelin G. Geredeli\thanks{%
email address: pguvengeredeli2@unl.edu. The research of P.G. Geredeli was
partially supported by the NSF grant DMS-1616425 and an Edith T. Hitz
Fellowship.} \\
Department of Mathematics, University of Nebraska-Lincoln, 68588 USA, \\
Department of Mathematics, Hacettepe University, Ankara-Turkey}
\title{ Stability Analysis of Coupled Structural Acoustics PDE Models under
Thermal Effects and with no Additional Dissipation }
\maketitle

\begin{abstract}
In this study we consider a coupled system of partial differential equations
(PDE's) which describes a certain structural acoustics interaction. One
component of this PDE system is a wave equation, which serves to model the
interior acoustic wave medium within a given three dimensional chamber $%
\Omega $. This acoustic wave equation is coupled on a boundary interface ($%
\Gamma _{0}$) to a two dimensional system of thermoelasticity: this
thermoelastic PDE comprises a structural beam or plate equation, which
governs the vibrations of flexible wall portion $\Gamma _{0}$ of the chamber 
$\Omega $; the elastic dynamics is coupled to a heat equation which also
evolves on $\Gamma _{0}$, and which imparts a thermal damping onto the
entire structural acoustic system. As we said, the interaction between the
wave and thermoelastic PDE components takes place on the boundary interface $%
\Gamma _{0}$, and involves coupling boundary terms which are above the level
of finite energy. We analyze the stability properties of this coupled
structural acoustics PDE model, in the absence of \ any additive feedback
dissipation on the hard walls $\Gamma _{1}$ of the boundary $\partial \Omega 
$. Under a certain geometric assumption on $\Gamma _{1}$, an assumption
which has appeared in the literature in conection with structural acoustic
flow, and which allows for the invocation of a recently derived microlocal
boundary trace estimate, we show that classical solutions of this thermally
damped structural acoustics PDE decay uniformly to zero, with a rational
rate of decay.

\bigskip

\end{abstract}


\section{Introduction}

\medskip

\subsection{Opening Remark}

In this work, we will consider stability properties of a coupled partial
differential equation (PDE) model which describes structural acoustic flow
under the influence of thermal dissipation: In particular, the PDE model
under consideration mathematically describes the interaction between an
interior acoustic field which evolves within a three dimensional chamber $%
\Omega $, and the structural displacements which occur along the flexible
portion $\Gamma _{0}$ of the chamber walls $\partial \Omega $, with these
elastic displacements being subjected to thermal effects. Consequently, the
PDE model consitutes a coupling of a wave equation and a two dimensional
system of thermoelasticity; the coupling between the two distinct dynamics
occurs on a boundary interface $\Gamma _{0}$. It is wellknown that
(uncoupled) thermoelastic plate systems, under all possible mechanical
boundary conditions, and with or without an accounting of rotational
inertia, have solutions which decay exponentially in time; see \cite{Liu}, 
\cite{Liu-2}, \cite{avalos4}, \cite{avalos5}, \cite{L-T.3}, \cite{L-T.4}. In
this connection, we are presently interested in discerning the extent to
which the (boundary) temperature dissipation in the thermoelastic PDE
component propagates onto the entire structural acoustic system,
particularly the interior acoustic wave component. In particular, our
objective is to study stability properties of said structural acoustic
systems, subject to thermal effects on boundary portion $\Gamma _{0}$, and
with \emph{no} additional feedback dissipation imposed on the
\textquotedblleft inactive\textquotedblright\ portion of the boundary,
denoted throughout as $\Gamma _{1}$.

\subsection{PDE Model}

Let $\Omega \subset \mathbb{R}^{3}$ be a bounded and open set with $C^{2}$%
-boundary $\partial \Omega =\Gamma =\overline{\Gamma }_{0}\cup \overline{%
\Gamma }_{1},$ where each $\Gamma _{i}$ is nonvoid, and $\Gamma _{0}\cap
\Gamma _{1}=\emptyset .$ In addition, the boundary segment $\Gamma _{0}$ is
flat. (In the statement of our main stability result, there will be two
additional assumptions made on geometry $\Omega $.) The wave equation is
invoked here to describe the interior acoustic medium within the three
dimensional chamber (or spatial domain) $\Omega $; this PDE dynamics is
coupled to a thermoelastic PDE which evolves on flat segment $\Gamma _{0}$.
To wit, the PDE system under consideration is given below, in solution
variables $\left[ z,z_{t},w,w_{t},\theta \right] $:%
\begin{align}
& \left\{ 
\begin{array}{l}
z_{tt}=\Delta z-z\text{ \ \ in \ }(0,T)\times 
\Omega
\\ 
\\ 
\frac{\partial z}{\partial \nu }=0\text{ \ \ on \ }(0,T)\times \Gamma _{1}
\\ 
\frac{\partial z}{\partial \nu }=w_{t}\text{ \ \ on \ }(0,T)\times \Gamma
_{0} \\ 
\\ 
\left[ z(0),z_{t}(0)\right] =\left[ z_{0},z_{1}\right] \in H^{1}(%
\Omega
)\times L^{2}(%
\Omega
)%
\end{array}%
\right.  \label{sys_w} \\
&  \notag \\
& \left\{ 
\begin{array}{l}
w_{tt}-\gamma \Delta w_{tt}+\Delta ^{2}w+\alpha \Delta \theta
+z_{t}|_{\Gamma _{0}}=0\text{ \ \ on \ }(0,T)\times \Gamma _{0} \\ 
\theta _{t}-\Delta \theta +\theta -\alpha \Delta w_{t}=0\text{ \ on \ }%
(0,T)\times \Gamma _{0} \\ 
\\ 
w=0\text{ \ \ on \ }(0,T)\times \partial \Gamma _{0} \\ 
\Delta w+(1-%
\mu
)B_{1}w+\alpha \theta =0\text{ \ \ on \ }(0,T)\times \partial \Gamma _{0} \\ 
\frac{\partial \theta }{\partial \nu }+\lambda \theta =0,\text{ }(\lambda
\geq 0)\text{ \ \ on \ }(0,T)\times \partial \Gamma _{0} \\ 
\\ 
\left[ w(0),w_{t}(0),\theta (0)\right] =\left[ w_{0},w_{1},\theta _{0}\right]
\in \lbrack H^{2}(\Gamma _{0})\cap H_{0}^{1}(\Gamma _{0})]\times H_{0,\gamma
}^{1}(\Gamma _{0})\times L^{2}(\Gamma _{0})%
\end{array}%
\right.  \label{sys_p}
\end{align}%
As given, $z$ is the wave solution component of the structural acoustic
model (\ref{sys_w})-(\ref{sys_p}), with quantity $z_{t}$ essentially
manifesting the underlying acoustic pressure of the chamber medium. The
dependent variable $w$ solves the elastic equation, with simply supported
boundary conditions, and evolves along the wall portion $\Gamma _{0}$. The
coupling between the two dynamics is accomplished via the Neumann boundary
condition for the $z$-wave equation, and via the Dirichlet trace $\left.
z_{t}\right\vert _{\Gamma \,_{0}}$, as it appears as a forcing term in the
Euler-Bernoulli (if $\gamma =0$) or Kirchoff plate (if $\gamma >0$) elastic
equation. In addition, the heat equation, in solution variable $\theta $,
with under either Robin or Neumann boundary conditions, manifests the
dissipation to which the structural acoustic system is subjected.

With regard to the physical parameters in the system of thermoelasticity in (%
\ref{sys_p}): the parameter $\gamma $ accounts for rotational forces, is
proportional to the square of the thickness of the plate, and is assumed to
be small, with $0\leq \gamma \leq 1$. In addition, the \textquotedblleft
coefficient of thermal expansion\textquotedblright\ $\alpha >0$ (see \cite%
{lagnese2}). The finite energy space component $H_{0,\gamma }^{1}(\Gamma
_{0})$, from which mechanical velocity data $w_{1}$\ is drawn in (\ref{sys_p}%
), and which depends on the value of $\gamma \geq 0$, is given by%
\begin{equation}
H_{0,\gamma }^{1}(\Gamma _{0})=\left\{ 
\begin{array}{l}
H_{0}^{1}(\Gamma _{0}),\text{ \ if }\gamma >0, \\ 
L^{2}(\Gamma _{0}\text{\ }),\text{ \ if }\gamma =0.%
\end{array}%
\right.  \label{drawn}
\end{equation}
Moreover, in connection with the mechanical boundary condition, the boundary
operator $B_{1}$ is given by%
\begin{equation*}
B_{1}w=2\nu _{1}\nu _{2}\frac{\partial ^{2}w}{\partial x\partial y}-\nu
_{1}^{2}\frac{\partial ^{2}w}{\partial y^{2}}-\nu _{2}^{2}\frac{\partial
^{2}w}{\partial x^{2}}
\end{equation*}%
where $\nu =(\nu _{1},\nu _{2})$ is the outer unit normal to the boundary,
and constant $0<%
\mu
<1$ is Poisson's modulus.

\subsection{Notation}

In this paper for a given domain $D$, its associated $L^{2}(D)$ will be
denoted as $||\cdot ||_{D}$. Inner products in $L^{2}(D)$ will be denoted by 
$\langle \cdot ,\cdot \rangle _{D}$. The space $H^{s}(D)$ will denote the
Sobolev space of order $s$, defined on a domain $D$, and $H_{0}^{s}(D)$
denotes the closure of $C_{0}^{\infty }(D)$ in the $H^{s}(D)$-norm $\Vert
\cdot \Vert _{H^{s}(D)}$ or $\Vert \cdot \Vert _{s,D}$.

\section{Abstract Setup}

The coupled system (\ref{sys_w})-(\ref{sys_p}) can be associated with a $%
C_{0}$-contraction semigroup on the space of initial data

\begin{equation}
\mathbf{H}=H^{1}(%
\Omega
)\times L^{2}(%
\Omega
)\times H_{0}^{2}(\Gamma _{0})\times H_{0,\gamma }^{1}(\Gamma _{0})\times
L^{2}(\Gamma _{0}),  \label{H}
\end{equation}%
with $\mathbf{H}$-inner product given by 
\begin{equation}
\begin{array}{l}
\left( \left[ z_{1},z_{2},w_{1},w_{2},\theta _{0}\right] ,\left[ \tilde{z}%
_{1},\tilde{z}_{2},\tilde{w}_{1},\tilde{w}_{2},\tilde{\theta}_{0}\right]
\right) _{\mathbf{H}}= \\ 
\text{ \ \ \ \ \ }\left( \nabla z_{1},\nabla \tilde{z}_{1}\right) _{\Omega
}+\left( z_{1},\tilde{z}_{1}\right) _{\Omega }+\left( z_{2},\tilde{z}%
_{2}\right) _{\Omega }+a\left( w_{1},\tilde{w}_{1}\right) _{\Gamma _{0}} \\ 
\text{ \ \ \ \ \ \ }+\left( w_{2},\tilde{w}_{2}\right) _{\Gamma _{0}}+\gamma
\left( \nabla w_{2},\nabla \tilde{w}_{2}\right) _{\Gamma _{0}}+\left( \theta
,\tilde{\theta}\right) _{\Gamma _{0}}.%
\end{array}
\label{H-I}
\end{equation}%
Here, bilinear form $a(\cdot ,\cdot ):H^{2}(\Gamma _{0})\times H^{2}(\Gamma
_{0})\rightarrow \mathbb{C}$ is that associated with the biharmonic
plate-bending problem; namely,%
\begin{equation}
a(\phi ,\psi )=\int_{\Gamma _{0}}\left[ \phi _{xx}\overline{\psi }_{xx}+\phi
_{yy}\overline{\psi }_{yy}+\mu (\phi _{xx}\overline{\psi }_{yy}+\phi _{yy}%
\overline{\psi }_{xx}+2(1-\mu )\phi _{xy}\overline{\psi }_{xy}\right]
d\Gamma _{0}  \label{plate}
\end{equation}%
(see \cite{lagnese2}).

\medskip

To explicitly describe this semigroup, and subsequently analyze its
stability properties, we will need to define the following operators, with
which we will abstractly model the PDE system (\ref{sys_w})-(\ref{sys_p}):

\begin{itemize}
\item Let the positive, self-adjoint operator $A_{N}:D(A_{N})\subset
L^{2}(\Omega )\rightarrow L^{2}(\Omega )$ be defined by

\begin{equation}
A_{N}\equiv I-\Delta ,\text{ \ \ \ }D(A_{N})=\{z\in H^{2}(\Omega ):\frac{%
\partial z}{\partial \nu }=0\text{ \ on \ }\Gamma \}.  \label{A_N}
\end{equation}

\item Associated with $A_{N}$ is the following harmonic extension of
boundary data, $N:L^{2}(\Gamma _{0})\rightarrow L^{2}(\Omega )$: 
\begin{equation}
Ng=h\Leftrightarrow \left\{ 
\begin{array}{c}
(I-\Delta )(h)=0\text{ \ in \ }\Omega , \\ 
\dfrac{\partial h}{\partial \nu }=\left\{ 
\begin{array}{c}
0\text{ \ on \ }\Gamma _{1}, \\ 
g\text{ \ on \ }\Gamma _{0}.\text{\ }%
\end{array}%
\right.%
\end{array}%
\right.  \label{N}
\end{equation}%
By elliptic regularity, for every $s\geq -1/2,$ $N\in \mathcal{L(}%
H^{s}(\Gamma _{0}),H^{s+\frac{3}{2}}(\Omega )).$

\item In addition, we define $\mathring{A}:D(\mathring{A})\subset
L^{2}(\Gamma _{0})\longrightarrow L^{2}(\Gamma _{0})$ by 
\begin{equation}
\mathring{A}f\equiv \Delta ^{2}f\text{, \ with }D(\mathring{A})=\{w\in
H^{4}(\Gamma _{0})\cap H_{0}^{1}(\Gamma _{0}):\Delta w+(1-%
\mu
)B_{1}w=0\text{ \ on \ }\partial \Gamma _{0}\}.  \label{angs}
\end{equation}

$\mathring{A}$ is a positive definite, self-adjoint operator, whose domain
of definition allows for the characterization 
\begin{equation}
D(\mathring{A}^{1/2})=H^{2}(\Gamma _{0})\cap H_{0}^{1}(\Gamma _{0})
\label{char_1}
\end{equation}%
(see \cite{grisvard}).

\item We also define $A_{D}:D(A_{D})\subset L^{2}(\Gamma _{0}$\ $%
)\longrightarrow L^{2}(\Gamma _{0}$\ $)$ to be the Laplacian operator with
homogeneous Dirichlet boundary conditions; namely, 
\begin{equation}
A_{D}=-\Delta ,\text{ }D(A_{D})=H^{2}(\Gamma _{0})\cap H_{0}^{1}(\Gamma
_{0});  \label{A_D}
\end{equation}%
By \cite{grisvard} we have the characterization%
\begin{equation}
D(A_{D})=H_{0}^{1}(\Gamma _{0}).  \label{char_2}
\end{equation}

\item Associated with $A_{D}$ is the operator $P_{\gamma }:$ $D(P_{\gamma
})\subset L^{2}(\Gamma _{0})\rightarrow L^{2}(\Gamma _{0})$, defined by%
\begin{equation}
P_{\gamma }=\left\{ 
\begin{array}{l}
I+\gamma A_{D},\text{ \ if }\gamma >0 \\ 
I,\text{ \ if }\gamma =0\text{,}%
\end{array}%
\right. \text{ \ \ with }D(P_{\gamma })=H^{2}(\Gamma _{0})\cap
H_{0}^{1}(\Gamma _{0}).  \label{P_g}
\end{equation}

Denoting $H_{0,\gamma }^{^{-1}}(\Gamma _{0})$ to be the topological dual of $%
H_{0,\gamma }^{^{1}}(\Gamma _{0})$, as defined in (\ref{drawn}), we then have%
\begin{equation}
H_{0,\gamma }^{1}(\Gamma _{0})=D(P_{\gamma }^{\frac{1}{2}})=\left\{ 
\begin{array}{c}
H_{0}^{1}(\Gamma _{0}),\text{ \ if }\gamma >0, \\ 
L^{2}(\Gamma _{0}),\text{ \ if }\gamma =0%
\end{array}%
\right. \text{ \ }\Rightarrow \text{ \ }P_{\gamma }\in \mathcal{L}\left(
H_{0,\gamma }^{1}(\Gamma _{0}),H_{0,\gamma }^{-1}(\Gamma _{0})\right) .
\label{H_P}
\end{equation}

\item In addition, $\gamma _{0}\in \mathcal{L}(H^{1}(\Gamma
_{0}),H^{1/2}(\partial \Gamma _{0}))$ will denote the classical Sobolev
trace map which yields for $f\in C^{\infty }(\overline{\Omega })$%
\begin{equation*}
\gamma _{0}(f)=f|_{\Gamma }.
\end{equation*}

\item We also have need of the following extensions of boundary data, which
are companion to $A_{D}$ and $\mathring{A}$, respectively:

(i) (Dirichlet) map $D:L^{2}(\partial \Gamma _{0})\rightarrow L^{2}(
\Gamma _{0})$ is given by 
\begin{equation}
Dh=v\Longleftrightarrow \left\{ 
\begin{array}{l}
\Delta v=0\text{ \ \ in \ }\Gamma _{0}, \\ 
v=h\text{ \ \ on \ }\partial \Gamma _{0},%
\end{array}%
\right.  \label{D}
\end{equation}

(ii) (Green) map $G:L^{2}(\partial \Gamma _{0})\longrightarrow L^{2}(\Gamma
_{0})$ is given by 
\begin{equation}
Gh=v\Longleftrightarrow \left\{ 
\begin{array}{l}
\Delta ^{2}v=0\text{ \ \ in \ }\Gamma _{0}, \\ 
v=0\text{ \ \ on \ }\partial \Gamma _{0}, \\ 
\Delta v+(1-%
\mu
)B_{1}v=h\text{ \ \ on \ }\partial \Gamma _{0}.%
\end{array}%
\right.  \label{G}
\end{equation}

By \cite{L-M}, we have for $s\in \mathbb{R}$, respectively, $D\in \mathcal{L}%
(H^{s}(\partial \Gamma _{0}),H^{s+\frac{1}{2}}(\Gamma _{0}))$ and $G\in 
\mathcal{L}(H^{s}(\partial \Gamma _{0}),H^{s+\frac{5}{2}}(\Gamma _{0})).$

\item Throughout, we will use repeatedly the fact that $N^{\ast }A_{N}\in 
\mathcal{L}(H^{1}(\Omega ),L^{2}(\Gamma _{0}))$ and \newline
$G^{\ast }\mathring{A}\in \mathcal{L}(H^{2}(\Gamma _{0})\cap
H_{0}^{1}(\Gamma _{0}),L^{2}(\Gamma _{0}))$ can be respectively
characterized as 
\begin{equation}
N^{\ast }A_{N}f=\left. f\right\vert _{\Gamma _{0}}\text{ \ for }f\in
D(A_{N}^{\frac{1}{2}})\text{, \ \ and \ \ }G^{\ast }\mathring{A}\varpi =%
\frac{\partial \varpi }{\partial n}\text{ \ for }\varpi \in D(\mathring{A}^{%
\frac{1}{2}})\text{,}  \label{char}
\end{equation}
\end{itemize}

(as can be verified outright by integrations by parts).

\bigskip


By means of the abstract operators defined above, the PDE model (\ref{sys_w}%
)-(\ref{sys_p}) with solution $[z,z_{t},w,w_{t},\theta ]$ can be concisely
rewritten as the following first order Cauchy problem:

\medskip

\begin{equation}
\frac{d}{dt}\left[ 
\begin{array}{c}
z \\ 
z_{t} \\ 
w \\ 
w_{t} \\ 
\theta%
\end{array}%
\right] =\mathbf{A}\left[ 
\begin{array}{c}
z \\ 
z_{t} \\ 
w \\ 
w_{t} \\ 
\theta%
\end{array}%
\right] \text{, \ \ }\left[ 
\begin{array}{c}
z(0) \\ 
z_{t}(0) \\ 
w(0) \\ 
w_{t}(0) \\ 
\theta (0)%
\end{array}%
\right] =\left[ 
\begin{array}{c}
z_{0} \\ 
z_{1} \\ 
w_{0} \\ 
w_{1} \\ 
\theta _{0}%
\end{array}%
\right] \in \mathbf{H}.  \label{ODE}
\end{equation}

\bigskip

\bigskip

Here the matrix operator $\mathbf{A}$ $:D(\mathbf{A})\subset \mathbf{H}%
\rightarrow \mathbf{H}$ is defined as%
\begin{equation}
\mathbf{A=}\left[ 
\begin{array}{ccccc}
0 & I & 0 & 0 & 0 \\ 
-A_{N} & 0 & 0 & A_{N}N & 0 \\ 
0 & 0 & 0 & I & 0 \\ 
0 & -P_{\gamma }^{-1}(\cdot )|_{\Gamma _{0}} & -P_{\gamma }^{-1}\mathring{A}
& 0 & \alpha P_{\gamma }^{-1}[A_{D}(I-D\gamma _{0})-\mathring{A}G\gamma _{0}]
\\ 
0 & 0 & 0 & -\alpha A_{D} & -A_{D}(I-D\gamma _{0})-I%
\end{array}%
\right]  \label{generator}
\end{equation}%
with 
\begin{equation}
\begin{array}{l}
D(\mathbf{A})=\{\left[ z_{1},z_{2},w_{1},w_{2},\theta _{0}\right] \in H^{1}(%
\Omega
)\times H^{1}(%
\Omega
)\times D(\mathring{A}^{1/2})\times D(\mathring{A}^{1/2})\times H^{2}(\Gamma
_{0}): \\ 
\text{ \ \ \ \ \textit{(D.i)}\ }z_{1}-Nw_{2}\in D(A_{N})\text{; \ \textit{%
(D.ii)} }\mathring{A}w_{1}+\alpha \mathring{A}G\gamma _{0}\theta_{0} \in
H_{0,\gamma }^{^{-1}}(\Gamma _{0})\text{; \ \textit{(D.iii)} }\frac{\partial
\theta _{0}}{\partial \nu }+\lambda \theta _{0}=0\text{ \ on \ }\partial
\Gamma _{0}\}.%
\end{array}
\label{domain}
\end{equation}

\bigskip

With a view of invoking the Lumer-Phillips Theorem - see e.g., p.14 of \cite%
{pazy} -- one can readily proceed to show that $\mathbf{A}$ $:D(\mathbf{A}%
)\subset \mathbf{H}\rightarrow \mathbf{H}$ is a maximal dissipative operator
(as has been done for other structural acoustic systems; see e.g.,\cite{jota}
and \cite{semi}). More precisely, we have the following statement of
existence and uniqueness for solutions for the the structural acoustics PDE
system (\ref{sys_w})-(\ref{sys_p}):

\begin{theorem}
\label{WP} The linear operator $\mathbf{A}$ $:D(\mathbf{A})\subset \mathbf{H}%
\rightarrow \mathbf{H}$, defined in (\ref{generator}), generates a $C_{0}$%
-semigroup $\left\{ e^{\mathbf{A}t}\right\} _{t\geq 0}$ of contractions on $%
\mathbf{H}$. Thus, if we denote initial data $\left[
z_{0},z_{1},w_{0},w_{1},\theta _{0}\right] $ and the solution $\left[
z,z_{t},w,w_{t},\theta \right] $ of (\ref{sys_w})-(\ref{sys_p}) (or
equivalently the ODE (\ref{ODE})) to be%
\begin{equation}
\begin{array}{l}
\Phi (\tau )=\left[ z(\tau ),z_{t}(\tau ),w(\tau ),w_{t}(\tau ),\theta (\tau
)\right] \text{ \ for all }\tau \geq 0 \\ 
\Phi _{0}=\Phi (0)=\left[ z_{0},z_{1},w_{0},w_{1},\theta _{0}\right] \text{,}%
\end{array}
\label{nomen}
\end{equation}%
then:

(a) $\Phi _{0}\in D(\mathbf{A})\Rightarrow \Phi (t)\in C([0,T];D(\mathbf{A}%
))\cap C^{1}([0,T];\mathbf{H})$ (continuously).

(b) $\Phi _{0}\in \mathbf{H}\Rightarrow \Phi (t)\in C([0,T];\mathbf{H})$
(continuously).

(c) Moreover, one has the following dissipative relation for all $t>0$:%
\begin{equation}
\int_{0}^{t}\left[ \left\Vert \nabla \theta \right\Vert _{\Gamma
_{0}}^{2}+\left\Vert \theta \right\Vert _{\Gamma _{0}}^{2}+\lambda
\left\Vert \theta \right\Vert _{\partial \Gamma _{0}}^{2}\right] d\tau
=\left\Vert \Phi _{0}\right\Vert _{\mathbf{H}}^{2}-\left\Vert \Phi
(t)\right\Vert _{\mathbf{H}}^{2}  \label{dissi}
\end{equation}%
(and so by the contraction of the semigroup $\left\{ e^{\mathbf{A}t}\right\}
_{t\geq 0}$, $\Phi _{0}\in \mathbf{H}\Rightarrow \theta \in L^{2}(0,\infty
;H^{1}(\Gamma _{0}))$, continuously).
\end{theorem}

\section{The Main Results}

\medskip

\subsection{Statement of Results}

\medskip

As a point of departure, we first note that finite energy solutions to the
system (\ref{sys_w})-(\ref{sys_p}) decay asymptotically in long time.
Indeed, the conclusion that the strong stability property holds for this
structural acoustics system can be made straightaway, having in hand the
dissipative relation (\ref{dissi}) and the compactness of the resolvent of $%
\mathbf{A}$ $:D(\mathbf{A})\subset \mathbf{H}\rightarrow \mathbf{H}$ (see
Corollary 3.1 of \cite{bench}; also \cite{slem} and \cite{levan}). And the
resolvent operator of $\mathbf{A}$ is indeed compact: In particular, we have
from (D.ii) of (\ref{domain}), the containment $D(\mathring{A}^{\frac{3}{4}%
})\subset H^{3}(\Gamma _{0})$ which is given in \cite{grisvard}, and the
regularity for elliptic operator $G$ as given in (\ref{G}), that%
\begin{equation*}
\left[ z_{1},z_{2},w_{1},w_{2},\theta _{0}\right] \in D(\mathbf{A}%
)\Longrightarrow w_{1}\in H^{3}(\Gamma _{0}).
\end{equation*}

The underlying thermal damping (\ref{dissi}), and the compactness of the
structural acoustic resolvent, give then the following decay result:

\begin{theorem}
\label{SS} With reference to the system (\ref{sys_w})-(\ref{sys_p}), the
continuous semigroup $\left\{ e^{\mathbf{A}t}\right\} _{t\geq 0}$ is
strongly stable. That is, with initial data $\Phi _{0}\in \mathbf{H}$, the
solution of (\ref{sys_w})-(\ref{sys_p}), or equivalently (\ref{ODE}), obeys $%
\lim_{t\rightarrow \infty }\left\Vert \Phi (t)\right\Vert _{\mathbf{H}}=0.$
\end{theorem}

\begin{remark}
In line with existing results in the literature, this strong stability is
obtained without any geometric assumptions on the boundary. The conclusion
of Theorem \ref{SS} can also be gleaned in what follows: to wit, in the
course of establishing our main result, Theorem \ref{US} below, it is
necessary to show that $i%
\mathbb{R}
\subset \rho (\mathbf{A})$. Consequently, Theorem \ref{SS} comes as a by
product, after invoking the wellknown spectral criteria for strong
asymptotic decay in \cite{A-B}.
\end{remark}

Our main result below will be valid under the following geometric
assumptions:

\medskip

\begin{description}
\item[\textsf{(Geometry.1)}] The inactive boundary portion $\Gamma _{1}$ is
convex;

\item[\textsf{(Geometry.2)}] There exists a point $x_{0}\in 
\mathbb{R}
^{3}$ such that $(x-x_{0})\cdot \nu \leq 0$ for all $x\in \Gamma _{1}.$
\end{description}



The bulk of our effort here is directed to proving the following result:
\begin{theorem}
\label{US} With $ \Omega \subset \mathbb{R}^{3} $ being a bounded and open set with $ C^2 $- boundary, and with no geometric assumptions, then $i%
\mathbb{R}
\subset \rho (\mathbf{A})$. Moreover, if assumptions \textsf{(Geometry.1)} and \textsf{%
(Geometry.2)} are in place, then there exist a positive constant $%
\mathfrak{B}$ and $C^{\ast }>0$ such that, for all $\beta \in 
\mathbb{R}
$ which satisfies $\left\vert \beta \right\vert \geq \mathfrak{B}$, the
resolvent operator $\mathcal{R}(i\beta ;\mathbf{A})=(i\beta \mathbf{I-A)}%
^{-1}$ obeys the estimate,%
\begin{equation}
\left\Vert \mathcal{R}(i\beta ;\mathbf{A})\Phi ^{\ast }\right\Vert _{\mathbf{%
H}}\leq C^{\ast }\left\vert \beta \right\vert ^{8}.  \label{41}
\end{equation}
\end{theorem}
Accepting for the time being the validity of Theorem \ref{US}, we have
immediately our stated objective:

\begin{corollary}
\label{main}Using the denotations in (\ref{nomen}), if $\Phi (t)\in
C([0,T];D(\mathbf{A}))\cap C^{1}([0,T];\mathbf{H})$ is a solution of (\ref%
{sys_w})-(\ref{sys_p}), corresponding to intial data $\Phi _{0}\in D(\mathbf{%
A})$, then one has the uniform polynomial decay estimate,%
\begin{equation*}
\left\Vert \Phi (t)\right\Vert _{\mathbf{H}}\leq \frac{M}{t^{\frac{1}{8}}}%
\left\Vert \Phi _{0}\right\Vert _{D(\mathbf{A})}\text{, for }t>0\text{, and
some }M>0.
\end{equation*}
\end{corollary}
\textbf{Proof of Corollary \ref{main}: }This comes from combining Theorem %
\ref{US} with the resolvent criterion for rational decay in \cite{tomilov}:

\smallskip

\begin{theorem}
\label{bt} Let $\left\{ T(t)\right\} _{t\geq 0}$ be a bounded $C_{0}-$%
semigroup on a Hilbert space $H$ with generator $A$ such that $i%
\mathbb{R}
\subset \rho (A).$ Then, for fixed $\alpha >0,$ the following are equivalent:%
\newline
(i) $\left\Vert \mathcal{R}(is;A)\right\Vert =\mathcal{O}(\left\vert
s\right\vert ^{\alpha }),$ $\left\vert s\right\vert \rightarrow \infty $; 

\smallskip 

\noindent (ii) $\left\Vert T(t)A^{-1}\right\Vert =\mathcal{O}(t^{-\frac{1}{%
\alpha }}\ ),$ $\ t\rightarrow \infty .$
\end{theorem}


\begin{remark}
The geometric conditions \textsf{(Geometry.1)} and \textsf{(Geometry.2)}
have been invoked before, in the context of controlling and stabilizing
structural acoustic flows; see \cite{jpma} and \cite{ag}; these assumptions
essentially dictate that the \textquotedblleft roof\textquotedblright\ $%
\Gamma _{1}$of the acoustic chamber $\Omega $ not be \textquotedblleft too
deep\textquotedblright . Since $\Gamma _{0}$ is only a portion of the
boundary, the imposition of geometric conditions on uncontrolled boundary
portion $\Gamma _{1}$ is fully expected; see \cite{bardos}. The assumptions 
\textsf{(Geometry.1)} and \textsf{(Geometry.2)} allow for the construction
of a special vector field -- this vector is generated outright in Appendix
II of \cite{LL}; see also \cite{LTZ} -- which, in tandem with the classic
wave equation identities, give rise to the static wave estimate Theorem \ref%
{start} below (see \cite{ag}).
\end{remark}


\begin{remark}
To provide as much clarity in our stability proof as circumstances will
allow, we consider the structural acoustics system (\ref{sys_w})-(\ref{sys_p}%
) with structural displacement $w(t)$ satisfying the simply supported
boundary conditions. However, our results here would also pertain to all
possible structural boundary conditions, including the case where $w$
satisfies the so-called free boundary conditions; i.e.,%
\begin{eqnarray*}
\Delta w+(1-%
\mu
)B_{1}w+\alpha \theta  &=&0\ \ \text{on}\ (0,T)\times \partial \Gamma _{0};
\\
\frac{\partial \Delta w}{\partial \nu }+(1-\mu )B_{2}w-\gamma \frac{\partial
w_{tt}}{\partial \nu }+\alpha \frac{\partial \theta }{\partial \nu } &=&0\ \ 
\text{on}\ (0,T)\times \partial \Gamma _{0},
\end{eqnarray*}%
where 
\begin{equation*}
B_{2}\phi =\frac{\partial }{\partial \tau }\left[ (\nu _{1}^{2}-\nu _{2}^{2})%
\frac{\partial ^{2}\phi }{\partial x_{1}\partial x_{2}}+\nu _{1}\nu _{2}(%
\frac{\partial ^{2}\phi }{\partial x_{2}^{2}}-\frac{\partial ^{2}\phi }{%
\partial x_{1}^{2}})\right] .
\end{equation*}%
(These mechanical boundary conditions were imposed in \cite{LL}). The
strategy employed here for the simply supported case would be successful for
the free case, although the details of proof for the latter would be more
cumbersome -- in particular, there would be a need to impose an additional
energy method so as to control the term $\left. \Delta w\right\vert
_{\partial \Gamma _{0}}$, as was undertaken in \cite{avalos4} for uncoupled
thermoelastic systems (see also \cite{lions} and \cite{L-T.2}, where such
mechanical boundary traces were first derived for plate dynamics).
\end{remark}


\subsection{Relevant Literature and Further Remarks}

\medskip

In this paper, our main goal is to ascertain uniform stability properties of
the coupled system (\ref{sys_w})-(\ref{sys_p}) for \emph{classical }%
solutions; i.e., for initial data $\Phi _{0}\in D(\mathbf{A})$; this is
Corollary \ref{main}. The important feature of this model is that it is
under the influence of thermal effects \emph{alone}, with no additive
feedback on the \textquotedblleft hard walls\textquotedblright\ $\Gamma _{1}$%
. This is in contrast to much of the existing literature in which uniform
stabilization results of various structural acoustic systems, for given
initial data of finite energy, have been attained by imposing additional
feedback damping on $\Gamma _{0}$, along with the intrinsic damping
mechanism coming from the elastic PDE component on $\Gamma _{0}$; see e.g., 
\cite{expo},\cite{fahroo},\cite{LL},\cite{grobb}. 

\smallskip 

In the literature, there do appear uniform decay results for structural
acoustics dynamics which do not have additional boundary dissipation on $%
\partial \Omega $. The earlier results required special initial data: For
example, \cite{litt} considered the spectral and stability properties of a
canonical wave -- damped second order ODE interaction on a rectangular
domain. The canonical setting in \cite{litt} allows for a Fourier analysis
which culminates in a statement of uniform decays for solutions which
correspond to certain smooth initial data (smoother than the domain of the
associated structural acoustic generator; see Theorem 5 of \cite{litt}). In
addition, the paper \cite{rebarber} gives an unspecified, uniform rate of
decay for solutions of a wellknown PDE model for structural acoustic flows,
for zero wave initial data. More recently, in \cite{ag}, rates of \emph{%
rational decay} were obtained for classical solutions of the aforesaid
wellknown structural acoustic PDE model: In \cite{ag}, the interior wave
equation in (\ref{sys_w}) appears as it does here -- the acoustic PDE
component does not seem to change from model to model (dating from \cite%
{beale}) -- however, the thermoelastic system (\ref{sys_p}) is replaced by
either a wave or plate equation under some degree of structural damping
(from [weak] viscous to [strong] Kelvin-Voight). In \cite{ag}, the primary
issue is appropriately dealing with the wave solution component for general
structural acoustic systems, as well as the boundary traces of both wave
displacement and velocity -- the main wave estimates of \cite{ag} are
applied in the present work. Subsequently, since the structural velocity in 
\cite{ag} manifests the damping to the entire PDE system, once we are able
to control the wave energy (in part by a $\Psi DO$ analysis of the wave PDE
component), much of the heavy lifting has been accomplished, since the
structural velocity dissipation directly contacts the wave component in \cite%
{ag}, via the wave Neumann boundary condition. The present situation stands
in stark contrast to that in \cite{ag}: the wave component of (\ref{sys_w})
is not directly coupled with the heat equation in (\ref{sys_p}), and so it
is not at all clear how (or if) the thermal damping alone in (\ref{sys_p})
will elicit some notion of uniform decay.

\smallskip 

The aforementioned paper \cite{LL} also considers the structural acoustic
system (\ref{sys_w})-(\ref{sys_p}), a PDE system under heat dissipation
(however with the mechanical component in \cite{LL} satisfying free boundary
conditions, rather than the simply supported which prevail here.) In \cite%
{LL}, the objective is to devise a (minimally invasive) dissipative feedback
scenario by which \emph{finite energy }solutions to a controlled version of (%
\ref{sys_w})-(\ref{sys_p}) decay uniformly, with respect to initial data in
the Hilbert space $\mathbf{H}$ of wellposedness. In this connection, this
earlier work imposes nonlinear boundary damping term $g(\left.
z_{t}\right\vert _{\Gamma _{1}})$ on hard walls $\Gamma _{1}$; however,
there is no such damping enforced on active boundary portion $\Gamma _{0}$.
In this way, \cite{LL} obtains explicit rates of decay; in particular, if
the nonlinearity can be \textquotedblleft bounded from below by a linear
function\textquotedblright , then the decay rate is of exponential type. In
contrast, our objective here is to investigate decay properties of solutions
to the structural acoustic model (\ref{sys_w})-(\ref{sys_p}), where the only
dissipation acting upon the system is that emanating from the thermal
gradient of the thermoelastic component on $\Gamma _{0}$. Given this
localized and indirect form of the damping -- in particular, the wave
component $[z,z_{t}]$ of (\ref{sys_w})-(\ref{sys_p}) is only influenced
indirectly by the temperature variable $\theta $ -- it would seem that one
should look to derive \emph{rational }rates of decay (as we did in \cite{ag}%
). Accordingly, our main result Theorem \ref{US} (and its Corollary \ref%
{main}) deal with the specification of polynomial rates of decay, which are
uniform with respect to smooth initial data; i.e., initial data which is
drawn from the domain of the associated structural acoustic generator $%
\mathbf{A}$ $:D(\mathbf{A})\subset \mathbf{H}\rightarrow \mathbf{H}$.

\smallskip

Structural acoustic PDE systems were initially considered by mathematical
control theorists with a view towards optimization with respect to the
implementation of (open loop) piezoelectric point controllers see e.g., \cite%
{banks} and \cite{jota}; stability properties of uncontrolled structural
acoustic dynamics were an afterthought. When the problem of stability for
solutions of structural acoustic systems -- with no additional dissipative
feedback -- was eventually taken up, it was quickly realized that such
stability analysis is generally a difficult problem (even the
\textquotedblleft soft\textquotedblright\ notion of strong decay is fairly
nontrivial in the case the associated structural acoustic generator has 
\emph{noncompact} resolvent; see \cite{semi} and \cite{litt}). Within the
context of the structural acoustic system (\ref{sys_w})-(\ref{sys_p}), we
can precisely state one serious complication: Using the notation (\ref{nomen}%
), if $\Phi (t)$ is a solution of (\ref{sys_w})-(\ref{sys_p}), with respect
to \emph{finite energy} data $\Phi _{0}\in \mathbf{H}$, then the wave trace
terms $(z_{t}|_{\Gamma _{0}})$ and $(\frac{\partial z}{\partial \tau }%
|_{\Gamma _{0}})$ are \textit{a priori} ill-defined in $L^{2}$-sense.
Consequently, it is extremely problematic to carry out the requisite time
domain multiplier methods , by way of generating the Lyapunov-type
inequalities which characterize uniform decay of finite energy solutions.
(See e.g., \cite{chen},\cite{lagnese},\cite{trigg} for instructive
illustrations in the context of the wave equation.) This finite energy
complication persists in the present case: although our main result Theorem %
\ref{US} gives a polynomial decay rate for \emph{smooth} solutions, the very
same wave boundary trace issues necessarily appear here, albeit it in the 
\emph{frequency}, and not \emph{time}, domain. This is owing to the fact
that we will work here to obtain the resolvent estimate in (\ref{41}) -- an
estimate in the \emph{finite energy }norm -- in accordance with the
resolvent criterion for rational decay of classical solutions to dissipative
systems, which is given in Theorem \ref{bt}. (This is also why a fair-sized
chunk of the paper is devoted to establishing that $i\mathbb{R}\subset \rho (%
\mathbf{A)}$.) However, a transformation of the system (\ref{sys_w})-(\ref%
{sys_p}) into the frequency domain does allow for the following: By way of
providing appropriate $L^{2}$-estimates for the wave boundary trace terms \ $%
(z_{t}|_{\Gamma _{0}})$ and $(\frac{\partial z}{\partial \tau }|_{\Gamma
_{0}})$ -- actually, to be precise, estimates for \textquotedblleft Laplace
transformed\textquotedblright\ versions of these wave boundary traces -- we
appeal in the course of our stability proof to the $\Psi DO$ result Theorem %
\ref{tang} below, which as we said was recently derived in \cite{ag} for
solutions of static wave equations, in the course of establishing rational
decay rates for solutions of a more canonical structural acoustic model
under the effects of structural damping. (To PDE control theorists, this
Theorem \ref{tang} should readily be recognized as a \textquotedblleft
frequency domain version\textquotedblright\ of the estimate in \cite{AMO},
which has been constantly used in the time domain for control of tangential
derivatives of wave equations.) We venture to say that if one undertakes the
stability analysis, \emph{in the frequency domain}, of any given
wave-structure PDE interaction, an he or she might find Theorem \ref{tang}
to be useful in some way; see e.g., \cite{avalos3},\cite{grobb}.

\smallskip

As we said, the geometric conditions \textsf{(Geometry.1)} and \textsf{%
(Geometry.2)} allow for the invocation of Theorem \ref{start}, which was
derived in \cite{ag} for the control of the wave component of a given
(static) structural acoustic system, regardless of the particular makeup of
the structural component. So the real issues in the present work, as it
concerns rational decay for smooth solutions of (\ref{sys_w})-(\ref{sys_p})
are to: (i) \ establish the containment $i\mathbb{R}\subset \rho (\mathbf{A)}
$; (ii) generate the needed energy inequalities for the structural component 
$[w,w_{t}]$, by way of ultimately obtaining the resolvent estimate (\ref{41}%
); (iii) completely understand how the thermal dissipation on $\Gamma _{0}$
allows for $L^{2}$-control of $z_{t}|_{\Gamma _{0}}$, and so then control of
the wave energy component. To deal with these issues (ii) and (iii), in the
course of proof of Theorem \ref{US}, we invoke the wellknown thermoelastic
multiplier $A_{D}^{-1}\theta $, where $A_{D}:D(A_{D}):L^{2}(\Gamma _{0}$ $%
)\rightarrow L^{2}(\Gamma _{0}$ $)$ (actually, we invoke the [formal]
Laplace transform of this multiplier); see, \cite{avalos4}, \cite{avalos5},%
\cite{expo},\cite{LL}. 

\smallskip 

Ultimately, we find that smooth solutions of (\ref{sys_w})-(\ref{sys_p})
decay at the rate $\mathcal{O}(t^{-\frac{1}{8}})$. As final note, we mention
that this is the same rational decay rate which is obtained in \cite{grobb2}%
, a work which deals with the rational decay problem for a thermoelastic
system composed of a Mindlin-Timoshenko plate (MT) which is subjected to a
thermal damping. Although the MT-heat system in \cite{grobb2} is not at all
associated with structural acoustic dynamics, one confronts the same
fundamental issue as in the present paper: in \cite{grobb2}, one of the
three \textquotedblleft shear angles\textquotedblright\ of the hyperbolic MT
components does not at all contact the the thermal component of the
dynamics; yet, as in the present paper, this situation of indirect damping
still allows for some measure of stability for solutions of the entire \cite%
{grobb2}.


\section{Associated Spectral Analysis}

\label{strong}


We begin by characterizing the inherent dissipativity of the PDE system (\ref%
{sys_w})-(\ref{sys_p}) in the frequency domain. Throughout we will use the
denotations for respective pre-images and images of the structural acoustic
generator (\ref{generator})-(\ref{domain}):%
\begin{equation}
\Phi \equiv \lbrack z_{1},z_{2},w_{1},w_{2},\theta ]\in D(\mathbf{A})\text{,
\ and data }\Phi ^{\ast }\equiv \lbrack z_{1}^{\ast },z_{2}^{\ast
},w_{1}^{\ast },w_{2}^{\ast },\theta ^{\ast }]\in \mathbf{H}.  \label{den}
\end{equation}

\smallskip

\begin{proposition}
\label{thermal}For given $\Phi \in D(\mathbf{A})$, as in (\ref{den}), we
have the relation%
\begin{equation}
Re\left( \mathbf{A}\Phi ,\Phi \right) _{\mathbf{H}}=-\left\Vert
\nabla \theta \right\Vert _{\Gamma _{0}}^{2}-\left\Vert \theta \right\Vert
_{\Gamma _{0}}^{2}-\lambda \left\Vert \theta \right\Vert _{\partial \Gamma
_{0}}^{2}.  \label{dissi_2}
\end{equation}
\end{proposition}
\textbf{Proof of Proposition \ref{thermal}:} Using the definition of the
inner product in (\ref{H-I}), and the explicit form of the matrix in (\ref%
{generator}), we have%
\begin{eqnarray}
\left( \mathbf{A}\Phi ,\Phi \right) _{\mathbf{H}} &=&\left( A_{N}^{\frac{1}{2%
}}z_{2},A_{N}^{\frac{1}{2}}z_{1}\right) _{\Omega }-\left\langle
A_{N}z_{1},z_{2}\right\rangle +\left\langle A_{N}Nw_{2},z_{2}\right\rangle 
\notag \\
&&+\left( \mathring{A}^{\frac{1}{2}}w_{2},\mathring{A}^{\frac{1}{2}%
}w_{1}\right) _{\Gamma _{0}}-\left( \left. z_{2}\right\vert _{\Gamma
_{0}},w_{2}\right) _{\Gamma _{0}}-\left\langle \mathring{A}%
w_{1},w_{2}\right\rangle  \notag \\
&&-\alpha \left( \Delta \theta ,w_{2}\right) _{\Gamma _{0}}-\alpha
\left\langle \mathring{A}G\gamma _{0}\theta ,w_{2}\right\rangle  \notag \\
&&+\alpha \left( \Delta w_{2},\theta \right) _{\Gamma _{0}}+\left( [\Delta
-I]\theta ,\theta \right) _{\Gamma _{0}}.  \label{mid}
\end{eqnarray}%
Upon further integrations by parts, and invocations of (\ref{char}), we have
now%
\begin{eqnarray*}
\left( \mathbf{A}\Phi ,\Phi \right) _{\mathbf{H}} &=&-2iIm\left(
A_{N}^{\frac{1}{2}}z_{1},A_{N}^{\frac{1}{2}}z_{2}\right) _{\Omega }-2i
Im\left( A_{N}^{\frac{1}{2}}w_{1},A_{N}^{\frac{1}{2}}w_{2}\right) _{\Omega }
\\
&&-2iIm\left( \left. z_{2}\right\vert _{\Gamma _{0}},w_{2}\right)
_{\Gamma _{0}}-2i\alpha Im\left( \gamma _{0}\theta ,\frac{\partial
w_{2}}{\partial n}\right) _{\partial \Gamma _{0}}-2iIm\left( \Delta
\theta ,w_{2}\right) _{\Gamma _{0}} \\
&&-\left\Vert \nabla \theta \right\Vert _{\Gamma _{0}}^{2}-\left\Vert \theta
\right\Vert _{\Gamma _{0}}^{2}-\lambda \left\Vert \theta \right\Vert
_{\partial \Gamma _{0}}^{2}.
\end{eqnarray*}%
Taking the real parts of both sides of this relation then establishes the
result. \ \ \ $\square $

\bigskip

In what follows, we will consider the abstract equation%
\begin{equation}
(i\beta I-\mathbf{A})\Phi =\Phi ^{\ast }\text{ \ for }\beta \in \mathbb{R},
\label{relation}
\end{equation}%
where $\Phi $ and $\Phi ^{\ast }$ are as given in (\ref{den}).

\bigskip

This section is devoted to proving the following requisite result on the
spectrum of the structural acoustic generator which is associated with the
system (\ref{sys_w})-(\ref{sys_p}), by way of ultimately establishing
Theorem \ref{US}.

\medskip

\begin{theorem}
\label{AB}With $\mathbf{A}$ $:D(\mathbf{A})\subset \mathbf{H}\rightarrow 
\mathbf{H}$ as defined in (\ref{generator})-(\ref{domain}), one has $i%
\mathbb{R}\subseteq \rho (\mathbf{A}).$
\end{theorem}
\textbf{Proof of Theorem \ref{AB}:}

\medskip 

\emph{Step 1. (}$0\in \rho (\mathbf{A})$\emph{)}. We will directly show that 
$\mathbf{A}$ has bounded inverse. With $\beta =0$ in (\ref{relation}), we
consider the equation 
\begin{equation}
\mathbf{A}\Phi =\Phi ^{\ast },  \label{zero}
\end{equation}%
where $\Phi $ and $\Phi ^{\ast }$ are as given in (\ref{den}). Using the
definition of the matrix in (\ref{generator}), then have the following
relations:%
\begin{equation}
z_{2}=z_{1}^{\ast },  \label{2}
\end{equation}%
\begin{equation}
-A_{N}z_{1}+A_{N}Nw_{2}=z_{2}^{\ast },  \label{3}
\end{equation}%
\begin{equation}
w_{2}=w_{1}^{\ast },  \label{4}
\end{equation}%
\begin{equation}
-P_{\gamma }^{-1}(z_{2})|_{\Gamma _{0}}-P_{\gamma }^{-1}\mathring{A}%
w_{1}+\alpha \lbrack P_{\gamma }^{-1}A_{D}(I-D\gamma _{0})-P_{\gamma }^{-1}%
\mathring{A}G\gamma _{0}]\theta =w_{2}^{\ast },  \label{5}
\end{equation}%
\begin{equation}
-\alpha A_{D}w_{2}-[A_{D}(I-D\gamma _{0})+I]\theta =\theta ^{\ast }.
\label{6}
\end{equation}%
Using (\ref{3}) and (\ref{4}), it is easy to see that 
\begin{equation}
z_{1}=Nw_{1}^{\ast }-A_{N}^{-1}z_{2}^{\ast }.  \label{8}
\end{equation}%
Also by (\ref{6}) and (\ref{4}), we have 
\begin{equation}
\theta =-\alpha A_{R}^{-1}A_{D}w_{1}^{\ast }-A_{R}^{-1}\theta ^{\ast },
\label{7}
\end{equation}%
where the positive definite, self-adjoint operator $A_{R}:D(A_{R})\subset
L^{2}(\Gamma _{0})\rightarrow $ $L^{2}(\Gamma _{0})$ is given by 
\begin{equation}
A_{R}=(I-\Delta )\text{, \ with }D(A_{R})=\{\theta \in H^{2}(\Gamma _{0}):%
\frac{\partial \theta }{\partial \nu }+\lambda \theta =0\text{ \ on \ }%
\partial \Gamma _{0}\}.  \label{A_R}
\end{equation}%
Lastly, using (\ref{2}) and (\ref{7}) in (\ref{5}) we take%
\begin{equation*}
w_{1}=-\mathring{A}^{-1}(z_{1}^{\ast })|_{\Gamma _{0}}+\alpha \lbrack 
\mathring{A}^{-1}A_{D}(I-D\gamma _{0})-G\gamma _{0}](-\alpha
A_{R}^{-1}A_{D}w_{1}^{\ast }-A_{R}^{-1}\theta ^{\ast })-\mathring{A}%
^{-1}P_{\gamma }w_{2}^{\ast }.
\end{equation*}%
Applying this relation, together with (\ref{2}), (\ref{4}), (\ref{8}) and (%
\ref{7}), we have that inverse operator $\mathbf{A}^{-1}\in \mathcal{L}(%
\mathbf{H},D(\mathbf{A}))$ exists, and is given by%
\begin{equation*}
\mathbf{A}^{-1}=\left[ 
\begin{array}{ccccc}
0 & -A_{N}^{-1} & N & 0 & 0 \\ 
I & 0 & 0 & 0 & 0 \\ 
-\mathring{A}^{-1}(\cdot )|_{\Gamma _{0}} & 0 & \Psi _{1} & -\mathring{A}%
^{-1}P_{\gamma } & \Psi _{2} \\ 
0 & 0 & I & 0 & 0 \\ 
0 & 0 & -\alpha A_{R}^{-1}A_{D} & 0 & -A_{R}^{-1}%
\end{array}%
\right] ,
\end{equation*}%
where 
\begin{equation*}
\Psi _{1}=-\alpha ^{2}[\mathring{A}^{-1}A_{D}(I-D\gamma _{0})-G\gamma
_{0}]A_{R}^{-1}A_{D},
\end{equation*}%
\begin{equation*}
\Psi _{2}=-\alpha \lbrack \mathring{A}^{-1}A_{D}(I-D\gamma _{0})-G\gamma
_{0}]A_{R}^{-1}.
\end{equation*}

\medskip

\emph{Step 2. (}$i\beta \notin \sigma _{p}(\mathbf{A})$, for $\beta \in 
\mathbb{R}
\diagdown \{0\}$\emph{).} If $\beta \neq 0$ and $\Phi ^{\ast }=\mathbf{0}$
in (\ref{relation}), then corresponding vector $\Phi \in D(\mathbf{A})$
satisfies 
\begin{equation}
\mathbf{A}\Phi =i\beta \Phi .  \label{eigen}
\end{equation}%
From the definition (\ref{generator}) this gives the the following equations:%
\begin{equation}
z_{2}=i\beta z_{1},  \label{9}
\end{equation}%
\begin{equation}
-A_{N}z_{1}+A_{N}Nw_{2}=i\beta z_{2},  \label{10}
\end{equation}%
\begin{equation}
w_{2}=i\beta w_{1},  \label{11}
\end{equation}%
\begin{equation}
-P_{\gamma }^{-1}(z_{2})|_{\Gamma _{0}}-P_{\gamma }^{-1}\mathring{A}%
w_{1}+\alpha P_{\gamma }^{-1}[A_{D}(I-D\gamma _{0})-\mathring{A}G\gamma
_{0}]\theta =i\beta w_{2},  \label{12}
\end{equation}%
\begin{equation}
-\alpha A_{D}w_{2}-[A_{D}(I-D\gamma _{0})+I]\theta =i\beta \theta .
\label{13}
\end{equation}

\medskip

Therewith: upon taking the $\mathbf{H}$-inner product of both sides of (\ref%
{eigen}), and invoking the relation in Proposition \ref{thermal}, we infer
that 
\begin{equation}
\theta =0.  \label{13.1}
\end{equation}%
In turn, applying this to the heat equation in (\ref{13}), we obtain%
\begin{equation}
w_{2}=0.  \label{13.2}
\end{equation}
This relation and (\ref{11}) give in turn,%
\begin{equation}
w_{1}=0.  \label{23}
\end{equation}%
Subsequently, applying (\ref{13.2}) and (\ref{23}) into (\ref{12}) gives%
\begin{equation}
z_{2}|_{\Gamma _{0}}=0.  \label{24}
\end{equation}

\bigskip

Since (\ref{24}) and (\ref{9}) yield that also $z_{1}|_{\Gamma _{0}}=0,$ we
obtain from this relation, and (\ref{9}), (\ref{10}) and (\ref{13.2}), that
variable $z_{1}$ satisfies the following overdetermined elliptic eigenvalue
problem: 
\begin{equation*}
\left\{ 
\begin{array}{c}
-A_{N}z_{1}=-\beta ^{2}z_{1},\text{ \ in \ }\Omega \\ 
z_{1}|_{\Gamma _{0}}=0,\text{ \ on \ }\Gamma _{0} \\ 
\frac{\partial z_{1}}{\partial \nu }=0\text{ \ on \ }\Gamma .%
\end{array}%
\right.
\end{equation*}%
In consequence, Holmgren's theorem gives that%
\begin{equation}
z_{1}=0,  \label{25}
\end{equation}%
which together with (\ref{9}) also yields that 
\begin{equation}
z_{2}=0.  \label{26}
\end{equation}%
As a result, combining (\ref{13.1}), (\ref{13.2}), (\ref{23}), (\ref{25})
and (\ref{26}) gives that solution $\Phi $ of (\ref{eigen}) is trivial.

\bigskip

\emph{Step 3 (}$i\beta \notin \sigma _{r}(\mathbf{A})$, for $\beta \in 
\mathbb{R}
\diagdown \{0\}$\emph{).} To show that $\lambda =i\beta $ is not in the
residual spectrum of \newline
$\mathbf{A}$ $:D(\mathbf{A})\subset \mathbf{H}\rightarrow \mathbf{H}$, where 
$\beta $ is nonzero, it is enough to show that $\lambda $ is not an
eigenvalue of its adjoint operator $\mathbf{A}^{\ast }:D(\mathbf{A}^{\ast
})\subset \mathbf{H}\rightarrow \mathbf{H}$ (see e.g., p. 127 of \cite%
{friedman}). This operator can be readily computed to be \ 
\begin{equation*}
\mathbf{A}^{\ast }\mathbf{=}\left[ 
\begin{array}{ccccc}
0 & -I & 0 & 0 & 0 \\ 
A_{N} & 0 & 0 & -A_{N}N & 0 \\ 
0 & 0 & 0 & -I & 0 \\ 
0 & P_{\gamma }^{-1}(\cdot )|_{\Gamma _{0}} & P_{\gamma }^{-1}\mathring{A} & 
0 & -\alpha \lbrack P_{\gamma }^{-1}A_{D}(I-D\gamma _{0})-P_{\gamma }^{-1}%
\mathring{A}G\gamma _{0}] \\ 
0 & 0 & 0 & \alpha A_{D} & -A_{D}(I-D\gamma _{0})-I%
\end{array}%
\right] ,
\end{equation*}%
with $D(\mathbf{A}^{\ast })=D(\mathbf{A})$, as given in (\ref{domain}).
Therewith, one can proceed to show, just as in the proof of Proposition \ref%
{thermal}, that for given $\Phi =[z_{1},z_{2},w_{1},w_{2},\theta ]\in D(%
\mathbf{A}^{\ast })$, we have the relation%
\begin{equation}
Re\left( \mathbf{A}^{\ast }\Phi ,\Phi \right) _{\mathbf{H}%
}=-\left\Vert \nabla \theta \right\Vert _{\Gamma _{0}}^{2}-\left\Vert \theta
\right\Vert _{\Gamma _{0}}^{2}-\lambda \left\Vert \theta \right\Vert
_{\partial \Gamma _{0}}^{2}.
\end{equation}%
With this dissipation in hand, we can proceed as in \emph{Step 2 }to show
that for $\beta \neq 0$, $i\beta \notin \sigma _{p}(\mathbf{A}^{\ast })$.

\medskip

\emph{Step 4 (}$i\beta \notin \sigma _{c}(\mathbf{A})$, for $\beta \in 
\mathbb{R}
\diagdown \{0\}$\emph{)}. By Theorem 2.27, p. 128 of \cite{friedman}, it is
enough to prove that $i\beta $ is not in the \emph{approximate spectrum} of $%
\mathbf{A}$. Suppose otherwise: then there exists a sequence of vectors $%
\Phi _{n}=[z_{1,n},z_{2,n},w_{1,n},w_{2,n},\theta _{n}]$ such that for every 
$n$,%
\begin{equation}
\left\Vert \Phi _{n}\right\Vert _{\mathbf{H}}=1\text{ }\forall \text{ }n%
\text{,\ \ \ \ \ and \ \ \ \ \ \ }\left\Vert (i\beta \mathbf{I-A)}\Phi
_{n}\right\Vert _{\mathbf{H}}<\frac{1}{n}.  \label{as1}
\end{equation}

Denoting%
\begin{equation}
(i\beta \mathbf{I-A)}\Phi _{n}=\Phi _{n}^{\ast }=\left[ 
\begin{array}{c}
z_{1,n}^{\ast } \\ 
z_{2,n}^{\ast } \\ 
w_{1,n}^{\ast } \\ 
w_{2,n}^{\ast } \\ 
\theta _{n}^{\ast }%
\end{array}%
\right] \in \mathbf{H},  \label{as2}
\end{equation}%
we have then the following relations:%
\begin{equation}
i\beta z_{1,n}-z_{2,n}=z_{1,n}^{\ast },  \label{27}
\end{equation}%
\begin{equation}
i\beta z_{2,n}+A_{N}z_{1,n}-A_{N}Nw_{2,n}=z_{2,n}^{\ast },  \label{28}
\end{equation}%
\begin{equation}
i\beta w_{1,n}-w_{2,n}=w_{1,n}^{\ast },  \label{29}
\end{equation}%
\begin{equation}
i\beta P_{\gamma }w_{2,n}+z_{2,n}|_{\Gamma _{0}}+\mathring{A}w_{1,n}+\alpha
\Delta \theta _{n}+\alpha \mathring{A}G\gamma _{0}\theta _{n}=P_{\gamma
}w_{2,n}^{\ast },  \label{30}
\end{equation}%
\begin{equation}
i\beta \theta _{n}-\Delta \theta _{n}-\alpha \Delta w_{2,n}+\theta
_{n}=\theta _{n}^{\ast }.  \label{31}
\end{equation}
If we take the $\mathbf{H}$-inner product of both sides of (\ref{as2}), and
subsequently invoke Proposition \ref{thermal}, we obtain

\begin{equation*}
\left\Vert \nabla \theta _{n}\right\Vert _{\Gamma _{0}}^{2}+\left\Vert
\theta _{n}\right\Vert _{\Gamma _{0}}^{2}+\lambda \left\Vert \theta
_{n}\right\Vert _{\partial \Gamma _{0}}^{2}=Re\left( \Phi _{n}^{\ast
},\Phi _{n}\right) _{\mathbf{H}}.
\end{equation*}%
Applying Cauchy-Schwartz to right hand side of this relation, and
subsequently using (\ref{as1})-(\ref{as2}), we get%
\begin{equation}
\lim_{n\rightarrow \infty }\left\Vert \theta _{n}\right\Vert _{H^{1}(\Gamma
_{0})}=0.\text{ }  \label{32}
\end{equation}%
Thereafter, we multiply (\ref{31}) by $w_{2,n}$ so as to have%
\begin{equation*}
\alpha \left\Vert \nabla w_{2,n}\right\Vert _{\Gamma _{0}}^{2}=-i\beta
\left( \theta _{n},w_{2,n}\right) _{\Gamma _{0}}-\left\langle \nabla \theta
_{n},\nabla w_{2,n}\right\rangle _{\Gamma _{0}}-\left\langle \theta
_{n},w_{2,n}\right\rangle _{\Gamma _{0}}+\left\langle \theta _{n}^{\ast
},w_{2,n}\right\rangle _{\Gamma _{0}}.
\end{equation*}%
For the right hand side of this relation, we can apply (\ref{29}) and (\ref%
{as1})-(\ref{as2}) to infer 
\begin{equation}
\lim_{n\rightarrow \infty }\left\Vert P_{\gamma }^{\frac{1}{2}%
}w_{2,n}\right\Vert _{\Gamma _{0}}=0.\text{ \ \ \ }  \label{33}
\end{equation}%
Subsequently, multiplying (\ref{30}) by $w_{1,n}$, and invoking (\ref{char}%
), we get%
\begin{equation*}
\left\Vert \mathring{A}^{\frac{1}{2}}w_{1,n}\right\Vert _{\Gamma
_{0}}^{2}=-i\beta \left( P_{\gamma }^{\frac{1}{2}}w_{2,n},P_{\gamma }^{\frac{%
1}{2}}w_{1,n}\right) _{\Gamma _{0}}-\left( z_{2,n}|_{\Gamma
_{0}},w_{1,n}\right) _{\Gamma _{0}}
\end{equation*}%
\begin{equation}
+\alpha \left( \nabla \theta _{n},\nabla w_{1,n}\right) _{\Gamma
_{0}}-\alpha \left( \theta _{n},\frac{\partial w_{1,n}}{\partial \nu }%
\right) _{\partial \Gamma _{0}}+\left( P_{\gamma }^{\frac{1}{2}%
}w_{2,n}^{\ast },P_{\gamma }^{\frac{1}{2}}w_{1,n}\right) _{\Gamma _{0}}.
\label{34}
\end{equation}%
By way of estimating the second term of RHS of (\ref{34}), we use (\ref{27})
and (\ref{29}) to have%
\begin{equation*}
\left\vert \left( z_{2,n}|_{\Gamma _{0}},w_{1,n}\right) _{\Gamma
_{0}}\right\vert \leq \left\vert \left( i\beta z_{1,n}|_{\Gamma _{0}},-\frac{%
i}{\beta }w_{2,n}-\frac{i}{\beta }w_{1,n}^{\ast }\right) _{\Gamma
_{0}}\right\vert +\left\vert \left( z_{1,n}^{\ast }|_{\Gamma _{0}},\frac{i}{%
\beta }w_{2,n}+\frac{i}{\beta }w_{1,n}^{\ast }\right) _{\Gamma
_{0}}\right\vert .
\end{equation*}%
Taking into account (\ref{33}) and (\ref{as1})-(\ref{as2}) in the last
relation gives now%
\begin{equation}
\lim_{n\rightarrow \infty }\left\vert \left( z_{2,n}|_{\Gamma
_{0}},w_{1,n}\right) _{\Gamma _{0}}\right\vert =0.  \label{35}
\end{equation}
For the fourth term of RHS of (\ref{34}), we apply Cauchy-Schwartz and the
Sobolev Trace Theorem to have%
\begin{eqnarray*}
\left\vert \left( \theta _{n},\frac{\partial w_{1,n}}{\partial \nu }\right)
_{\partial \Gamma _{0}}\right\vert &\leq &C\left\Vert \theta _{n}\right\Vert
_{\partial \Gamma _{0}}\left\Vert \frac{\partial w_{1,n}}{\partial \nu }%
\right\Vert _{\partial \Gamma _{0}} \\
&\leq &C\left\Vert \theta _{n}\right\Vert _{H^{1}(\Gamma _{0})}\left\Vert
w_{1,n}\right\Vert _{H^{2}(\Gamma _{0})}.
\end{eqnarray*}%
After considering (\ref{as1}) and (\ref{32}) , we then infer%
\begin{equation}
\lim_{n\rightarrow \infty }\left\vert \left( \theta _{n},\frac{\partial
w_{1,n}}{\partial \nu }\right) _{\partial \Gamma _{0}}\right\vert =0.
\label{36}
\end{equation}
From the relation (\ref{34}), in combination with (\ref{32}), (\ref{33}), (%
\ref{35}), (\ref{36}) and (\ref{as1})-(\ref{as2}), we then obtain 
\begin{equation}
\lim_{n\rightarrow \infty }\left\Vert \mathring{A}^{\frac{1}{2}%
}w_{1,n}\right\Vert _{\Gamma _{0}}=0.  \label{37}
\end{equation}
If we subsequently read off the relation (\ref{30}), and use (\ref{27}), (%
\ref{32}),(\ref{33}), and (\ref{37}), and (\ref{as1})-(\ref{as2}), we then
have the convergence%
\begin{equation}
\lim_{n\rightarrow \infty }z_{1,n}|_{\Gamma _{0}}=0\text{ \ \ in \ }\left[ D(%
\mathring{A}^{\frac{1}{2}})\right] ^{\prime }.  \label{38}
\end{equation}

\medskip

We must now deal with the wave component of $\Phi _{n}$: Using resolvent
relations (\ref{27}) and (\ref{28}) we obtain that the sequence $\left\{
z_{1,n}\right\} $ satisfies%
\begin{equation}
-\beta ^{2}z_{1,n}+A_{N}z_{1,n}-A_{N}Nw_{2,n}=z_{2,n}^{\ast }+i\beta
z_{1,n}^{\ast }.  \label{39}
\end{equation}

\smallskip

Moreover, from (\ref{as1}) we know that some subsequence $\left\{
z_{1,n}\right\} $ converges weakly to $\mathfrak{z}$, say, in $H^{1}(\Omega
).$ With this weak convergence in mind, we take the inner product of both
sides of (\ref{39}) with respect to a given $\psi \in H^{1}(\Omega )$, and
so have%
\begin{equation*}
-\beta ^{2}\left( z_{1,n},\psi \right) _{\Omega }+\left( A_{N}^{\frac{1}{2}%
}z_{1,n},A_{N}^{\frac{1}{2}}\psi \right) _{\Omega }-\left( w_{2,n},\psi
|_{\Gamma _{0}}\right) _{\Gamma _{0}}=\left( z_{2,n}^{\ast }+i\beta
z_{1,n}^{\ast },\psi \right) _{\Omega }\text{, }
\end{equation*}%
where we used the characterization in (\ref{char}). Now if we pass to the
weak limit above, simultaneously using (\ref{as1})-(\ref{as2}) and (\ref{33}%
), we see that weak limit $\mathfrak{z}\in H^{1}(\Omega )$ satisfies 
\begin{equation}
\left\langle (1-\beta ^{2})\mathfrak{z},\psi \right\rangle _{\Omega
}+\left\langle \nabla \mathfrak{z},\nabla \psi \right\rangle _{\Omega }=0%
\text{, \ \ \ for every }\psi \in H^{1}(\Omega ).  \label{over1}
\end{equation}%
In addition, by Rellich-Kondrasov, and the boundedness of the Sobolev Trace
Map, and the convergence in (\ref{38}), we infer that weak limit $\mathfrak{z%
}\in H^{1}(\Omega )$ of $\left\{ z_{1,n}\right\} $ has zero boundary trace
on $\Gamma _{0}$. This fact and (\ref{over1}) means that $\mathfrak{z}$
satisfies the following overdetermined problem:%
\begin{equation*}
\left\{ 
\begin{array}{c}
(1-\beta ^{2})\mathfrak{z}-\Delta \mathfrak{z}=0\text{ \ \ in \ }\Omega , \\ 
\frac{\partial \mathfrak{z}}{\partial \nu }=0\text{ \ \ on \ }\Gamma , \\ 
\mathfrak{z}=0\text{ \ \ on \ }\Gamma _{0}.%
\end{array}%
\right. 
\end{equation*}%
A subsequent application of Holmgren's Theorem gives that necessarily $%
\mathfrak{z}\equiv 0$; and so%
\begin{equation*}
z_{1,n}\overset{weakly}{\rightarrow }\text{ }0\text{ \ \ in \ }H^{1}(\Omega
).\text{\ }
\end{equation*}%
In turn, by the resolvent relation (\ref{27}) and the Rellich-Kondrasov
Theorem, we have that

\begin{equation}
\lim_{n\rightarrow \infty }\left\Vert z_{2,n}\right\Vert _{\Omega }=0.\text{%
\ }  \label{40}
\end{equation}%
To finish \emph{Step 4}: we multiply (\ref{28}) by $z_{1,n}$ to get%
\begin{equation*}
\left\Vert A_{N}^{\frac{1}{2}}z_{1,n}\right\Vert _{\Gamma _{0}}^{2}=-i\beta
\left( z_{2,n},z_{1,n}\right) _{\Omega }+\left( w_{2,n},z_{1,n}|_{\Gamma
_{0}}\right) _{\Gamma _{0}}+\left( z_{2,n}^{\ast },z_{1,n}\right) _{\Omega }
\end{equation*}%
(and also using (\ref{char})). Passing to limit as $n\rightarrow \infty $
and using (\ref{as1})-(\ref{as2}), (\ref{33}),(\ref{40}) and the Sobolev
Trace Theorem, we obtain%
\begin{equation}
\lim_{n\rightarrow \infty }\left\Vert A_{N}^{\frac{1}{2}}z_{1,n}\right\Vert
_{\Gamma _{0}}=0.  \label{40.5}
\end{equation}
The limits (\ref{32}), (\ref{33}), (\ref{37}), (\ref{40}) and (\ref{40.5})
gives now that 
\begin{equation*}
\underset{n\rightarrow \infty }{\lim }\left\Vert \Phi _{n}\right\Vert _{H}=0.
\end{equation*}%
This limit contradicts (\ref{as1}), and so no nonzero parameter $i\beta $ is
in the approximate spectrum of \newline
$\mathbf{A}$ $:D(\mathbf{A})\subset \mathbf{H}\rightarrow \mathbf{H}$. The
proof of Theorem \ref{AB} is now complete. \ \ \ $\square $

\medskip 

\section{Proof of Theorem \protect\ref{US}}

In this section, we give the proof of the resolvent estimate (\ref{41}),
which characterizes rational decay of the given structural acoustic dynamics.

Having proved that $i%
\mathbb{R}
\subseteq \rho (\mathbf{A})$ in Theorem \ref{AB}, we will now look at the
action of the resolvent operator on the imaginary axis. To this end, we
consider the equation 
\begin{equation}
(i\beta \mathbf{I-A)}\Phi =\Phi ^{\ast },  \label{ResRel}
\end{equation}%
where pre-image $\Phi $ and image $\Phi ^{\ast }$ are as given in (\ref{den}%
). As such, we have then the following relations: 
\begin{equation}
i\beta z_{1}-z_{2}=z_{1}^{\ast },  \label{42}
\end{equation}%
\begin{equation}
i\beta z_{2}+A_{N}z_{1}-A_{N}Nw_{2}=z_{2}^{\ast },  \label{43}
\end{equation}%
\begin{equation}
i\beta w_{1}-w_{2}=w_{1}^{\ast },  \label{44}
\end{equation}%
\begin{equation}
i\beta P_{\gamma }w_{2}+z_{2}|_{\Gamma _{0}}+\mathring{A}w_{1}+\alpha \Delta
\theta +\alpha \mathring{A}G\gamma _{0}\theta =P_{\gamma }w_{2}^{\ast },
\label{45}
\end{equation}%
\begin{equation}
i\beta \theta -\Delta \theta -\alpha \Delta w_{2}+\theta =\theta ^{\ast }.
\label{46}
\end{equation}

\medskip

\emph{Step I. }We start by obtaining an estimate on the thermal component of
the solution: we take the $\mathbf{H}$-inner product of both sides of (\ref%
{ResRel}) with respect to $\Phi $, and subsequently invoke Proposition \ref%
{thermal}. This gives the relation 
\begin{equation}
\left\Vert \nabla \theta \right\Vert _{\Gamma _{0}}^{2}+\left\Vert \theta
\right\Vert _{\Gamma _{0}}^{2}+\lambda \left\Vert \theta \right\Vert
_{\partial \Gamma _{0}}^{2}=\left\vert Re\left( \Phi ^{\ast },\Phi
\right) _{\mathbf{H}}\right\vert .  \label{teta}
\end{equation}

\bigskip

\emph{Step II (A preliminary estimate for the wave component of }$\Phi $%
\emph{). }This step would really be invariant with respect to the interior
wave component of any structural acoustic system under analysis for
polynomial decay properties; see e.g., the models considered in the
stability papers \cite{grobb}, \cite{LL}, and \cite{ag}. Using resolvent
relations (\ref{42}) and (\ref{44}) in (\ref{43}), we obtain the following
boundary value problem in $z_{1}:$%
\begin{equation}
\left\{ 
\begin{array}{l}
-\beta ^{2}z_{1}-\Delta z_{1}+z_{1}=z_{2}^{\ast }+i\beta z_{1}^{\ast }\text{
\ \ in \ }\Omega , \\ 
\frac{\partial z_{1}}{\partial \nu }=0\text{ \ on }\Gamma _{1}, \\ 
\frac{\partial z_{1}}{\partial \nu }=i\beta w_{1}-w_{1}^{\ast }\text{ \ on }%
\Gamma _{0}%
\end{array}%
\right.  \label{47}
\end{equation}

\bigskip

With respect to $z_{1}$, we apply here a preliminary estimate for wave
components of static structural acoustic systems which was recently derived,
under said geometric assumptions \textsf{(Geometry.1)} and \textsf{%
(Geometry.2).}

\medskip

\begin{theorem}
\label{start}(\cite[See Lemma 8 and inequality (50) therein.]{ag})\ \ Let
the geometric assumptions \textsf{(Geometry.1)} and \textsf{(Geometry.2)} be
in place. Then the wave component $z_{1}$ of the resolvent relation (\ref%
{ResRel}) -- or what is the same, the static wave equation (\ref{47}) --
obeys the following estimate, for $\left\vert \beta \right\vert $
sufficiently large, and arbitrary $\epsilon ^{\ast }>0$:%
\begin{equation}
\int\limits_{\Omega }\left\vert \nabla z_{1}\right\vert ^{2}d\Omega +\beta
^{2}\int\limits_{\Omega }\left\vert z_{1}\right\vert ^{2}d\Omega \leq
C\left( \left\Vert \frac{\partial z_{1}}{\partial \tau }\right\Vert _{\Gamma
_{0}}^{2}+\beta ^{2}\left\Vert z_{1}\right\Vert _{\Gamma _{0}}^{2}\right)
+\epsilon ^{\ast }\left\Vert \Phi \right\Vert _{H}^{2}+C\beta ^{2}\left\Vert
\Phi ^{\ast }\right\Vert _{\mathbf{H}}^{2}.  \label{pre_1}
\end{equation}
\end{theorem}

\bigskip

In order to render the wave estimate (\ref{pre_1}) useful, we must control
the tangential derivative of $z_{1}$ on right hand side (a term which is
strictly above the $H^{1}$-energy level for the wave displacement). To this
end, we appeal to the recently derived $\Psi DO$-result in \cite[Theorem 9]%
{ag}, whose wellknown time varying progenitor is in \cite{AMO}:

\medskip

\begin{theorem}
\label{tang}(See \cite[Theorem 9]{ag}). Let $\ \Gamma _{\ast }$ be a smooth
connected subset of boundary $\Gamma .$ Then the structural wave component
of the resolvent relation (\ref{ResRel}) -- or what is the same, the static
wave equation (\ref{47}) -- obeys the following boundary estimate, for
arbitrary $\delta >0:$%
\begin{equation}
\left\Vert \frac{\partial z_{1}}{\partial \tau }\right\Vert _{\Gamma _{\ast
}}\leq C^{\ast }\left\{ \left\Vert \beta z_{1}\right\Vert _{\Gamma _{\ast
}}+\left\vert \beta \right\vert \left\Vert w_{1}\right\Vert _{\Gamma
_{0}}+\left\Vert z_{1}\right\Vert _{H^{\frac{1}{2}+\delta }(\Omega
)}+\left\vert \beta \right\vert \left\Vert \Phi ^{\ast }\right\Vert _{%
\mathbf{H}}\right\} .  \label{z1-tang}
\end{equation}
\end{theorem}

\medskip

Utilizing the tangential estimate (\ref{z1-tang}) with respect to the right
hand side of (\ref{pre_1}) (with $\Gamma _{\ast }=\Gamma _{0}$ therein), we
now get the following initial estimate for the wave displacement component
of $\Phi =[z_{1},z_{2},w_{1},w_{2},\theta ]$:

\medskip

\begin{equation}
\int\limits_{\Omega }\left\vert \nabla z_{1}\right\vert ^{2}d\Omega +\beta
^{2}\int\limits_{\Omega }\left\vert z_{1}\right\vert ^{2}d\Omega \leq
\epsilon ^{\ast }\left\Vert \Phi \right\Vert _{\mathbf{H}}^{2}+C\left[
\left\Vert \beta z_{1}\right\Vert _{\Gamma _{0}}^{2}+\left\vert \beta
\right\vert ^{2}\left\Vert P_{\gamma }^{\frac{1}{2}}w_{1}\right\Vert
_{\Gamma _{0}}^{2}+\left\Vert z_{1}\right\Vert _{H^{\frac{1}{2}+\delta
}(\Omega )}^{2}+\left\vert \beta \right\vert ^{2}\left\Vert \Phi ^{\ast
}\right\Vert _{\mathbf{H}}^{2}\right] .  \label{z1-est}
\end{equation}

\bigskip

\emph{Step III. (The thermoelastic component).}

\medskip

We start by incorporating the resolvent relations (\ref{42}) and (\ref{44})
into (\ref{45}) (while recalling the operators in (\ref{angs}), (\ref{A_D}),
(\ref{P_g}), and (\ref{G})): we have%
\begin{equation}
-\beta ^{2}(I-\gamma \Delta )w_{1}+i\beta z_{1}|_{\Gamma _{0}}+\Delta
^{2}w_{1}+\alpha \Delta \theta =\mathcal{F}_{\beta }^{\ast },  \label{48}
\end{equation}%
where%
\begin{equation}
\mathcal{F}_{\beta }^{\ast }=P_{\gamma }w_{2}^{\ast }+i\beta P_{\gamma
}w_{1}^{\ast }+z_{1}^{\ast }|_{\Gamma _{0}}.  \label{F*}
\end{equation}

\medskip

In addition, inserting (\ref{44}) into (\ref{46}), we obtain the heat
equation%
\begin{equation}
i\beta \theta -i\alpha \beta \Delta w_{1}-\Delta \theta +\theta =\mathcal{G}%
^{\ast },  \label{49}
\end{equation}%
where%
\begin{equation*}
\mathcal{G}^{\ast }=\theta ^{\ast }-\alpha \Delta w_{1}^{\ast }.
\end{equation*}

\medskip

The PDE's (\ref{48}) and (\ref{49}), together with their respective boundary
conditions, then give the following inhomogeneous thermoelastic system:%
\begin{eqnarray}
&&\left\{ 
\begin{array}{l}
-\beta ^{2}(I-\gamma \Delta )w_{1}+i\beta z_{1}|_{\Gamma _{0}}+\Delta
^{2}w_{1}+\alpha \Delta \theta =\mathcal{F}_{\beta }^{\ast }\text{ \ in }%
\Gamma _{0}\text{\ ,} \\ 
i\beta \theta -i\alpha \beta \Delta w_{1}-\Delta \theta +\theta =\mathcal{G}%
^{\ast }\text{ \ \ in }\Gamma _{0},%
\end{array}%
\right.  \label{50} \\
&&  \notag \\
&&\left\{ 
\begin{array}{l}
w_{1}=0\text{ \ \ on \ }\partial \Gamma _{0}, \\ 
\Delta w_{1}+(1-%
\mu
)B_{1}w_{1}+\alpha \theta =0\text{ \ \ on \ }\partial \Gamma _{0}, \\ 
\frac{\partial \theta }{\partial \nu }+\lambda \theta =0,\text{ }(\lambda
\geq 0)\text{ \ \ on \ }\partial \Gamma _{0}.%
\end{array}%
\right.  \label{51}
\end{eqnarray}

\medskip

To obtain an estimate for the structural velocity term $\beta P_{\gamma }^{%
\frac{1}{2}}w_{1}$, we apply the operator $-iA_{D}^{-1}$ to the heat
equation in (\ref{50}), so as to have%
\begin{equation}
\beta A_{D}^{-1}\theta =-\alpha \beta w_{1}+i(I-D\gamma _{0})\theta
+iA_{D}^{-1}\theta -iA_{D}^{-1}\mathcal{G}^{\ast }.  \label{53}
\end{equation}

\smallskip

Therewith, we have%
\begin{equation}
\left( -\beta ^{2}(I-\gamma \Delta )w_{1},A_{D}^{-1}\theta \right) _{\Gamma
_{0}}=\alpha \beta ^{2}\left\Vert P_{\gamma }^{\frac{1}{2}}w_{1}\right\Vert
_{\Gamma _{0}}^{2}+i\beta \left( P_{\gamma }w_{1},[I-D\gamma _{0}]\theta
+A_{D}^{-1}\theta -A_{D}^{-1}\mathcal{G}^{\ast }\right) _{\Gamma _{0}}.
\label{54}
\end{equation}

\medskip

The relation (\ref{54}) gives then the estimate%
\begin{eqnarray}
\alpha \beta ^{2}\left\Vert P_{\gamma }^{\frac{1}{2}}w_{1}\right\Vert
_{\Gamma _{0}}^{2} &\leq &\left\vert \beta \left( P_{\gamma
}w_{1},[I-D\gamma _{0}]\theta +A_{D}^{-1}\theta -A_{D}^{-1}\mathcal{G}^{\ast
}\right) _{\Gamma _{0}}\right\vert \text{ \ \ \ \ \ \ \ \ \ \ \ \ \ \ }%
\boxed{E_{1}}  \notag \\
&&+\left\vert \left( \beta ^{2}(I-\gamma \Delta )w_{1},A_{D}^{-1}\theta
\right) _{\Gamma _{0}}\right\vert \text{ \ \ \ \ \ \ \ \ \ \ \ \ \ \ \ \ }%
\boxed{E_{2}}  \label{55}
\end{eqnarray}

\medskip

\textit{For the expression} $E_{1}$: From the definition of operator $D$ in (%
\ref{D}), and the characterization (\ref{H_P}), we have that $[I-D\gamma
_{0}]\in \mathcal{L}(H^{1}(\Gamma _{0}),D(P_{\gamma }^{\frac{1}{2}}))$.
Thus, 
\begin{equation}
\left\vert \beta \left( P_{\gamma }w_{1},[I-D\gamma _{0}]\theta \right)
_{\Gamma _{0}}\right\vert \leq \left\vert \beta \right\vert \left\Vert
P_{\gamma }^{\frac{1}{2}}w_{1}\right\Vert _{\Gamma _{0}}\left\Vert \theta
\right\Vert _{H^{1}(\Gamma _{0})}.  \label{r1}
\end{equation}%
Moreover, from the definition of $P_{\gamma }$ in (\ref{P_g}), we have that $%
P_{\gamma }A_{D}^{-1}\in \mathcal{L}([D(P_{\gamma }^{\frac{1}{2}})]^{\prime
},D(P_{\gamma }^{\frac{1}{2}}))$. Combining this boundedness with the
estimate (\ref{r1}) and $\left\vert ab\right\vert \leq \epsilon
a^{2}+C_{\epsilon }b^{2}$, we then have 
\begin{equation}
E_{1}\leq \epsilon \beta ^{2}\left\Vert P_{\gamma }^{\frac{1}{2}%
}w_{1}\right\Vert _{\Gamma _{0}}^{2}+C_{\epsilon }\left[ \left\Vert \theta
\right\Vert _{H^{1}(\Gamma _{0})}^{2}+\left\Vert \Phi ^{\ast }\right\Vert
_{H}^{2}\right] ,  \label{e1}
\end{equation}%
where data $\Phi ^{\ast }=[z_{1}^{\ast },z_{2}^{\ast },w_{1}^{\ast
},w_{2}^{\ast },\theta ^{\ast }]$ as denoted in (\ref{den}).

\medskip

\textit{For the expression} $E_{2}$: A transposing of self-adjoint $%
P_{\gamma }^{\frac{1}{2}}$ gives%
\begin{equation*}
\left( \beta ^{2}(I-\gamma \Delta )w_{1},A_{D}^{-1}\theta \right) _{\Gamma
_{0}}=\beta ^{2}\left( P_{\gamma }^{\frac{1}{2}}w_{1},P_{\gamma }^{\frac{1}{2%
}}A_{D}^{-1}\theta \right) _{\Gamma _{0}},
\end{equation*}%
whence we obtain, as $A_{D}^{-1}\in \mathcal{L}(L^{2}(\Gamma
_{0}),H^{2}(\Gamma _{0})\cap H_{0}^{1}(\Gamma _{0}))$,

\begin{equation}
E_{2}\leq \epsilon \beta ^{2}\left\Vert P_{\gamma }^{\frac{1}{2}%
}w_{1}\right\Vert _{\Gamma _{0}}^{2}+C_{\epsilon }\beta ^{2}\left\Vert
\theta \right\Vert _{\Gamma _{0}}^{2}.  \label{e1.5}
\end{equation}

\bigskip

Applying (\ref{e1}) and (\ref{e1.5}) to (\ref{55}), and taking $0<\epsilon
<\alpha /4$, we get now%
\begin{equation}
\beta ^{2}\left\Vert P_{\gamma }^{\frac{1}{2}}w_{1}\right\Vert _{\Gamma
_{0}}^{2}\leq C_{\epsilon ,\alpha }\left( \beta ^{2}\left\Vert \theta
\right\Vert _{H^{1}(\Gamma _{0})}^{2}+\left\Vert \Phi ^{\ast }\right\Vert _{%
\mathbf{H}}^{2}\right) .  \label{e5}
\end{equation}

\bigskip

By way of handling the mechanical displacement, we first recall the Green's
formula in \cite{lagnese2}: namely, for functions $\phi $ and $\psi $
sufficiently smooth, we have%
\begin{equation}
\left( \Delta ^{2}\phi ,\psi \right) _{\Gamma _{0}}=a(\phi ,\psi
)+\int\limits_{\partial \Gamma _{0}}\left[ \frac{\partial \Delta \phi }{%
\partial \nu }+(1-\mu )B_{2}\phi \right] \overline{\psi }-\int\limits_{%
\partial \Gamma _{0}}\left[ \Delta \phi +(1-\mu )B_{1}\phi \right] \frac{%
\partial \bar{\psi}}{\partial \nu };  \label{green}
\end{equation}%
where bilinear form $a(\cdot ,\cdot ):H^{2}(\Gamma _{0})\times H^{2}(\Gamma
_{0})\rightarrow \mathbb{C}$ is as given in (\ref{plate}), and boundary
expression 
\begin{equation*}
B_{2}\phi =\frac{\partial }{\partial \tau }\left[ (\nu _{1}^{2}-\nu _{2}^{2})%
\frac{\partial ^{2}\phi }{\partial x_{1}\partial x_{2}}+\nu _{1}\nu _{2}(%
\frac{\partial ^{2}\phi }{\partial x_{2}^{2}}-\frac{\partial ^{2}\phi }{%
\partial x_{1}^{2}})\right] .
\end{equation*}

\bigskip

Therewith, we multiply both sides of the structural equation in (\ref{50})
by $w_{1}$, integrate, and subsequently apply the Green's formula (\ref%
{green}). This gives%
\begin{eqnarray}
&&\left\Vert \mathring{A}^{\frac{1}{2}}w_{1}\right\Vert _{\Gamma _{0}}^{2} 
\notag \\
&=&-i\beta \left( z_{1}|_{\Gamma _{0}},w_{1}\right) _{\Gamma _{0}}-\alpha
\left( \Delta \theta ,w_{1}\right) _{\Gamma _{0}}-\alpha \left( \theta ,%
\frac{\partial w_{1}}{\partial n}\right) _{\partial \Gamma _{0}}+\beta
^{2}\left( P_{\gamma }w_{1},w_{1}\right) _{\Gamma _{0}}+\left( \mathcal{F}%
_{\beta }^{\ast },w_{1}\right) _{\Gamma _{0}}  \notag \\
&=&-i\beta \left( z_{1}|_{\Gamma _{0}},w_{1}\right) _{\Gamma _{0}}+\alpha
\left( \nabla \theta ,\nabla w_{1}\right) _{\Gamma _{0}}-\alpha \left(
\theta ,\frac{\partial w_{1}}{\partial n}\right) _{\partial \Gamma
_{0}}+\beta ^{2}\left\Vert P_{\gamma }^{\frac{1}{2}}w_{1}\right\Vert
_{\Gamma _{0}}^{2}+\left( \mathcal{F}_{\beta }^{\ast },w_{1}\right) _{\Gamma
_{0}},  \label{57}
\end{eqnarray}

where in the last step we integrated by parts once more.

\smallskip

We focus on the first term on the right hand side of (\ref{57}): Using the
heat equation in (\ref{50}), we have%
\begin{equation}
\left( z_{1}|_{\Gamma _{0}},i\beta w_{1}\right) _{\Gamma _{0}}=-\frac{1}{%
\alpha }\left( z_{1}|_{\Gamma _{0}},A_{D}^{-1}[i\beta \theta +\theta -%
\mathcal{G}^{\ast }]\right) _{\Gamma _{0}}-\frac{1}{\alpha }\left(
z_{1}|_{\Gamma _{0}},[I-D\gamma _{0}]\theta \right) _{\Gamma _{0}}.
\label{57.1}
\end{equation}

\medskip

For the first term on right hand side of (\ref{57.1}), we re-use the
structural equation in (\ref{50}): For $\left\vert \beta \right\vert >1$ we
have%
\begin{eqnarray}
&&-\frac{1}{\alpha }\left( z_{1}|_{\Gamma _{0}},A_{D}^{-1}[i\beta \theta
+\theta -\mathcal{G}^{\ast }]\right) _{\Gamma _{0}}  \notag \\
&=&\frac{1}{\alpha }\left( \beta ^{2}P_{\gamma }w_{1}-\Delta
^{2}w_{1}-\alpha \Delta \theta +\mathcal{F}_{\beta }^{\ast
},A_{D}^{-1}[\theta -\frac{i}{\beta }\theta +\frac{i}{\beta }\mathcal{G}%
^{\ast }]\right) _{\Gamma _{0}}.  \label{57.15}
\end{eqnarray}

\smallskip

Estimating right hand side of (\ref{57.15}), by means of (\ref{e5}), (\ref%
{green}), the Sobolev Trace Theorem and Cauchy-Schwartz, we have for $%
\left\vert \beta \right\vert >1$,%
\begin{equation}
\left\vert \frac{1}{\alpha }\left( z_{1}|_{\Gamma _{0}},A_{D}^{-1}[i\beta
\theta +\theta -\mathcal{G}^{\ast }]\right) _{\Gamma _{0}}\right\vert \leq
C\left( \beta ^{2}\left\Vert \theta \right\Vert _{H^{1}(\Gamma
_{0})}^{2}+\left\Vert \Phi ^{\ast }\right\Vert _{\mathbf{H}}^{2}+\left\Vert
\Phi \right\Vert _{\mathbf{H}}\left[ \left\Vert \theta \right\Vert _{\Gamma
_{0}}+\left\Vert \Phi ^{\ast }\right\Vert _{\mathbf{H}}\right] \right) .
\label{57.2}
\end{equation}

\medskip

For the second term on the right hand side of (\ref{57.1}), we have by the
boundedness in (\ref{D}), Cauchy-Schwartz and the Sobolev Trace Theorem,%
\begin{equation}
\left\vert \frac{1}{\alpha }\left( z_{1}|_{\Gamma _{0}},[I-D\gamma
_{0}]\theta \right) _{\Gamma _{0}}\right\vert \leq C\left\Vert
z_{1}\right\Vert _{H^{\frac{1}{2}+\delta }(\Omega )}\left\Vert \theta
\right\Vert _{H^{1}(\Gamma _{0})}.  \label{57.3}
\end{equation}

\bigskip

Applying (\ref{57.2}) and (\ref{57.3}) to the right hand side of (\ref{57.1}%
), we thus have for $\left\vert \beta \right\vert >1$,%
\begin{equation}
\left\vert \left( z_{1}|_{\Gamma _{0}},i\beta w_{1}\right) _{\Gamma
_{0}}\right\vert \leq C\left( \left\Vert z_{1}\right\Vert _{H^{\frac{1}{2}%
+\delta }(\Omega )}\left\Vert \theta \right\Vert _{H^{1}(\Gamma
_{0})}+\left\Vert \Phi \right\Vert _{\mathbf{H}}\left[ \left\Vert \theta
\right\Vert _{H^{1}(\Gamma _{0})}+\left\Vert \Phi ^{\ast }\right\Vert _{%
\mathbf{H}}\right] +\beta ^{2}\left\Vert \theta \right\Vert _{H^{1}(\Gamma
_{0})}^{2}+\left\Vert \Phi ^{\ast }\right\Vert _{\mathbf{H}}^{2}\right) .
\label{57.4}
\end{equation}

\medskip

Using in turn this estimate to majorize the right hand side of (\ref{57}),
along with (\ref{e5}), Cauchy-Schwartz and the Sobolev Trace Theorem, we
have then for $\left\vert \beta \right\vert >1$,%
\begin{equation}
\left\Vert \mathring{A}^{\frac{1}{2}}w_{1}\right\Vert _{\Gamma _{0}}^{2}\leq
C\left( \left\Vert z_{1}\right\Vert _{H^{\frac{1}{2}+\delta }(\Omega
)}\left\Vert \theta \right\Vert _{H^{1}(\Gamma _{0})}+\left\Vert \Phi
\right\Vert _{\mathbf{H}}\left[ \left\Vert \theta \right\Vert _{H^{1}(\Gamma
_{0})}+\left\Vert \Phi ^{\ast }\right\Vert _{\mathbf{H}}\right] +\beta
^{2}\left\Vert \theta \right\Vert _{H^{1}(\Gamma _{0})}^{2}+\left\Vert \Phi
^{\ast }\right\Vert _{\mathbf{H}}^{2}\right) .  \label{57.5}
\end{equation}

\bigskip

\emph{Step IV. (An appropriate estimate for the boundary traces }$\beta
z_{1}|_{\Gamma _{0}}$\emph{).}

\medskip

Given the right hand side (\ref{z1-est}), we apparently need control of $%
\left\Vert \beta z_{1}|_{\Gamma _{0}}\right\Vert _{\Gamma _{0}}$; in turn,
this boundary estimate will allow us to ultimately refine the right hand
side of (\ref{z1-est}), and subsequently recover the required energy
estimate for all components of the solution to (\ref{ResRel}). To this end,
we start by reconsidering the structural equation in (\ref{50}): We have,
upon applying the inverse $-i\mathring{A}^{-\frac{1}{2}}$ to both sides,%
\begin{equation}
\mathring{A}^{-\frac{1}{2}}\beta z_{1}|_{\Gamma _{0}}=-i\beta ^{2}\mathring{A%
}^{-\frac{1}{2}}P_{\gamma }w_{1}+i\mathring{A}^{\frac{1}{2}}w_{1}+\alpha i%
\mathring{A}^{\frac{1}{2}}G\gamma _{0}\theta -i\alpha \mathring{A}^{-\frac{1%
}{2}}A_{D}(I-D\gamma _{0})\theta -i\mathring{A}^{-\frac{1}{2}}\mathcal{F}%
_{\beta }^{\ast }.  \label{ex}
\end{equation}

\smallskip

Since $[D(A_{D}^{1/2})]^{\prime }=H^{-1}(\Gamma _{0})=[D(\mathring{A}%
^{1/4})]^{\prime }$ from \cite{grisvard}, this characterization and the
boundedness posted in (\ref{D}) and (\ref{G}), respectively, give the
initial estimate, for $\left\vert \beta \right\vert >1$,%
\begin{equation*}
\left\Vert \beta z_{1}|_{\Gamma _{0}}\right\Vert _{[D(\mathring{A}%
^{1/2})]^{\prime }}\leq C\left( \beta ^{2}\left\Vert P_{\gamma
}^{1/2}w_{1}\right\Vert _{\Gamma _{0}}+\left\Vert \mathring{A}^{\frac{1}{2}%
}w_{1}\right\Vert _{\Gamma _{0}}+\left\Vert \theta \right\Vert
_{H^{1}(\Gamma _{0})}+\left\vert \beta \right\vert \left\Vert \Phi ^{\ast
}\right\Vert _{\mathbf{H}}\right)
\end{equation*}

\smallskip

Applying the estimates (\ref{e5}) and (\ref{57.5}) to right hand side gives
now, for $\left\vert \beta \right\vert >1$, \ 
\begin{equation}
\begin{array}{l}
\left\Vert \beta z_{1}|_{\Gamma _{0}}\right\Vert _{[D(\mathring{A}%
^{1/2})]^{\prime }} \\ 
\text{ \ \ }\leq C\left( \sqrt{\left\Vert z_{1}\right\Vert _{H^{\frac{1}{2}%
+\delta }(\Omega )}\left\Vert \theta \right\Vert _{H^{1}(\Gamma _{0})}}+%
\sqrt{\left\Vert \Phi \right\Vert _{\mathbf{H}}\left\Vert \theta \right\Vert
_{H^{1}(\Gamma _{0})}}+\sqrt{\left\Vert \Phi \right\Vert _{\mathbf{H}%
}\left\Vert \Phi ^{\ast }\right\Vert _{\mathbf{H}}}+\left\vert \beta
\right\vert \left\Vert \theta \right\Vert _{H^{1}(\Gamma _{0})}+\left\vert
\beta \right\vert \left\Vert \Phi ^{\ast }\right\Vert _{\mathbf{H}}\right) .%
\end{array}
\label{f11}
\end{equation}

\medskip

We use this estimate in an interpolation between $H^{-2}(\Gamma _{0})$ and $%
H^{1}(\Gamma _{0})$ (see Theorem 12.4, p. 73 of \cite{L-M}): Namely we have 
\begin{eqnarray*}
\left\Vert \beta z_{1}|_{\Gamma _{0}}\right\Vert _{\Gamma _{0}} &\leq
&C\left\Vert \beta z_{1}|_{\Gamma _{0}}\right\Vert _{H^{-2}(\Gamma _{0})}^{%
\frac{1}{3}}\left\Vert \beta z_{1}|_{\Gamma _{0}}\right\Vert _{H^{1}(\Gamma
_{0})}^{\frac{2}{3}} \\
&\leq &C\left\vert \beta \right\vert ^{\frac{2}{3}}\left\Vert \beta
z_{1}|_{\Gamma _{0}}\right\Vert _{H^{-2}(\Gamma _{0})}^{\frac{1}{3}%
}\left\Vert z_{1}|_{\Gamma _{0}}\right\Vert _{H^{1}(\Gamma _{0})}^{\frac{2}{3%
}} \\
&\leq &C_{\epsilon _{1}}\left\vert \beta \right\vert ^{2}\left\Vert \beta
z_{1}|_{\Gamma _{0}}\right\Vert _{H^{-2}(\Gamma _{0})}+\frac{\epsilon _{1}}{%
C^{\ast }}\left\Vert z_{1}|_{\Gamma _{0}}\right\Vert _{H^{1}(\Gamma _{0})},
\end{eqnarray*}%
after also Young's Inequality, where $C^{\ast }$ is the constant which
appears in the tangential estimate (\ref{z1-tang}). Invoking Theorem \ref%
{tang} - with $\Gamma _{\ast }=\Gamma _{0}$ therein - we further obtain, for 
$\left\vert \beta \right\vert >1$, 
\begin{eqnarray}
&&\left\Vert \beta z_{1}|_{\Gamma _{0}}\right\Vert _{\Gamma _{0}}  \notag \\
&\leq &C_{\epsilon _{1}}\left\vert \beta \right\vert ^{2}\left\Vert \beta
z_{1}|_{\Gamma _{0}}\right\Vert _{H^{-2}(\Gamma _{0})}+\epsilon
_{1}\left\Vert \beta z_{1}\right\Vert _{\Gamma _{0}}+C_{\gamma ,\epsilon
_{1}}\left( \left\Vert z_{1}\right\Vert _{H^{\frac{1}{2}+\delta }(\Omega
)}+\left\vert \beta \right\vert \left\Vert P_{\gamma }^{\frac{1}{2}%
}w_{1}\right\Vert _{H^{1}(\Gamma _{0})}+\left\vert \beta \right\vert
\left\Vert \Phi ^{\ast }\right\Vert _{\mathbf{H}}\right)  \notag \\
&\leq &C_{\epsilon _{1}}\left\vert \beta \right\vert ^{2}\left\Vert \beta
z_{1}|_{\Gamma _{0}}\right\Vert _{H^{-2}(\Gamma _{0})}+\epsilon
_{1}\left\Vert \beta z_{1}\right\Vert _{\Gamma _{0}}+C_{\gamma ,\epsilon
_{1}}\left( \left\Vert z_{1}\right\Vert _{H^{\frac{1}{2}+\delta }(\Omega
)}+\left\vert \beta \right\vert \left\Vert \theta \right\Vert _{H^{1}(\Gamma
_{0})}+\left\vert \beta \right\vert \left\Vert \Phi ^{\ast }\right\Vert _{%
\mathbf{H}}\right) ,  \label{f2}
\end{eqnarray}%
where in the last step we invoked estimate (\ref{e5}).

\bigskip

For the first term on right hand side of (\ref{f2}), we invoke (\ref{f11})
(and $\left\vert ab\right\vert \leq \epsilon ^{\ast }a^{2}+C_{\epsilon
^{\ast }}b^{2}$), so as to have for $\left\vert \beta \right\vert >1$, 
\begin{eqnarray}
&&C_{\epsilon _{1}}\left\vert \beta \right\vert ^{2}\left\Vert \beta
z_{1}|_{\Gamma _{0}}\right\Vert _{H^{-2}(\Gamma _{0})}  \notag \\
&\leq &C\left\vert \beta \right\vert ^{2}\left( \sqrt{\left\Vert
z_{1}\right\Vert _{H^{\frac{1}{2}+\delta }(\Omega )}\left\Vert \theta
\right\Vert _{H^{1}(\Gamma _{0})}}+\sqrt{\left\Vert \Phi \right\Vert _{%
\mathbf{H}}\left\Vert \theta \right\Vert _{H^{1}(\Gamma _{0})}}+\sqrt{%
\left\Vert \Phi \right\Vert _{\mathbf{H}}\left\Vert \Phi ^{\ast }\right\Vert
_{\mathbf{H}}}+\left\vert \beta \right\vert \left\Vert \theta \right\Vert
_{H^{1}(\Gamma _{0})}+\left\vert \beta \right\vert \left\Vert \Phi ^{\ast
}\right\Vert _{\mathbf{H}}\right)  \notag \\
&\leq &\frac{\epsilon ^{\ast }}{2}\left\Vert \Phi \right\Vert _{\mathbf{H}%
}+C_{\epsilon ^{\ast }}\left( \left\Vert z_{1}\right\Vert _{H^{\frac{1}{2}%
+\delta }(\Omega )}+\left\vert \beta \right\vert ^{4}\left\Vert \theta
\right\Vert _{H^{1}(\Gamma _{0})}+\left\vert \beta \right\vert
^{4}\left\Vert \Phi ^{\ast }\right\Vert _{\mathbf{H}}\right) .  \label{f3}
\end{eqnarray}

(This is the point in the proof where the decay rate of the structural
acoustics model (\ref{sys_w})-(\ref{sys_p}) is determined.)

\medskip

Applying (\ref{f3}) to the right hand side of (\ref{f2}), and taking $%
0<\epsilon _{1}<1/2$, we have then the following controlling estimate for $%
\beta z_{1}|_{\Gamma _{0}}$ in $L^{2}$-topology:%
\begin{equation}
\left\Vert \beta z_{1}|_{\Gamma _{0}}\right\Vert _{\Gamma _{0}}\leq \epsilon
^{\ast }\left\Vert \Phi \right\Vert _{\mathbf{H}}+C_{\epsilon ^{\ast
}}\left( \left\Vert z_{1}\right\Vert _{H^{\frac{1}{2}+\delta }(\Omega
)}+\left\vert \beta \right\vert ^{4}\left\Vert \theta \right\Vert
_{H^{1}(\Gamma _{0})}+\left\vert \beta \right\vert ^{4}\left\Vert \Phi
^{\ast }\right\Vert _{\mathbf{H}}\right) .  \label{bz}
\end{equation}

\bigskip

\emph{Step V (Conclusion of the proof of Theorem \ref{US}). }Applying (\ref%
{bz}) to the right hand side of (\ref{z1-est}), re-invoking (\ref{e5}), and
rescaling $\epsilon ^{\ast }>0$, we have for $\left\vert \beta \right\vert
>1 $,%
\begin{equation}
\int\limits_{\Omega }\left\vert \nabla z_{1}\right\vert ^{2}d\Omega +\beta
^{2}\int\limits_{\Omega }\left\vert z_{1}\right\vert ^{2}d\Omega \leq
\epsilon ^{\ast }\left\Vert \Phi \right\Vert _{\mathbf{H}}^{2}+\widetilde{C}%
\left[ \left\Vert z_{1}\right\Vert _{H^{\frac{1}{2}+\delta }(\Omega
)}^{2}+\left\vert \beta \right\vert ^{8}\left\Vert \theta \right\Vert
_{H^{1}(\Gamma _{0})}^{2}+\left\vert \beta \right\vert ^{8}\left\Vert \Phi
^{\ast }\right\Vert _{\mathbf{H}}^{2}\right] .  \label{plu}
\end{equation}

\medskip

To deal with the lower order wave term on right hand side, we interpolate:
For $\left\vert \beta \right\vert >1$ we have 
\begin{eqnarray}
\left\Vert z_{1}\right\Vert _{H^{\frac{1}{2}+\delta }(\Omega )} &\leq
&\left\vert \frac{\beta }{\beta }\right\vert ^{\frac{1}{2}-\delta
}\left\Vert z_{1}\right\Vert _{\Omega }^{\frac{1}{2}-\delta }\left\Vert
z_{1}\right\Vert _{H^{1}(\Omega )}^{\frac{1}{2}+\delta }  \notag \\
&\leq &\sqrt{\frac{\epsilon _{2}}{2\widetilde{C}}}\left\Vert \beta
z_{1}\right\Vert _{\Omega }+C_{\epsilon _{2}}\frac{1}{\left\vert \beta
\right\vert ^{\frac{1-2\delta }{1+2\delta }}}\left\Vert z_{1}\right\Vert
_{H^{1}(\Omega )},  \label{per}
\end{eqnarray}%
after using Young's Inequality, where positive constant $\widetilde{C}$ is
the constant in (\ref{plu}). Applying this estimate to (\ref{plu}), we then
have for $\left\vert \beta \right\vert >1$,%
\begin{eqnarray*}
&&\left( \int\limits_{\Omega }\left\vert \nabla z_{1}\right\vert ^{2}d\Omega
+\int\limits_{\Omega }\left\vert z_{1}\right\vert ^{2}d\Omega \right)
+(\beta ^{2}-1)\int\limits_{\Omega }\left\vert z_{1}\right\vert ^{2}d\Omega
\\
&\leq &\epsilon ^{\ast }\left\Vert \Phi \right\Vert _{\mathbf{H}%
}^{2}+\epsilon _{2}\left\Vert \beta z_{1}\right\Vert _{\Omega
}^{2}+C_{\epsilon _{2}}\left[ \frac{1}{\left\vert \beta \right\vert ^{\frac{%
2-4\delta }{1+2\delta }}}\left\Vert z_{1}\right\Vert _{H^{1}(\Omega
)}^{2}+\left\vert \beta \right\vert ^{8}\left\Vert \theta \right\Vert
_{H^{1}(\Gamma _{0})}^{2}+\left\vert \beta \right\vert ^{8}\left\Vert \Phi
^{\ast }\right\Vert _{\mathbf{H}}^{2}\right] .
\end{eqnarray*}

\medskip

Taking $0<\epsilon _{2}<\frac{1}{2}$, we then have for $\left\vert \beta
\right\vert >1$,%
\begin{eqnarray*}
&&\left( \int\limits_{\Omega }\left\vert \nabla z_{1}\right\vert ^{2}d\Omega
+\int\limits_{\Omega }\left\vert z_{1}\right\vert ^{2}d\Omega \right)
+\left( \frac{1}{2}-\epsilon _{2}\right) \beta ^{2}\int\limits_{\Omega
}\left\vert z_{1}\right\vert ^{2}d\Omega \\
&& \\
&\leq &\epsilon ^{\ast }\left\Vert \Phi \right\Vert _{\mathbf{H}%
}^{2}+C_{\epsilon _{2}}\left[ \frac{1}{\left\vert \beta \right\vert ^{\frac{%
2-4\delta }{1+2\delta }}}\left\Vert z_{1}\right\Vert _{H^{1}(\Omega
)}^{2}+\left\vert \beta \right\vert ^{8}\left\Vert \theta \right\Vert
_{H^{1}(\Gamma _{0})}^{2}+\left\vert \beta \right\vert ^{8}\left\Vert \Phi
^{\ast }\right\Vert _{\mathbf{H}}^{2}\right] ,
\end{eqnarray*}%
whence we obtain for $\left\vert \beta \right\vert >1$,%
\begin{eqnarray*}
&&\left( \int\limits_{\Omega }\left\vert \nabla z_{1}\right\vert ^{2}d\Omega
+\int\limits_{\Omega }\left\vert z_{1}\right\vert ^{2}d\Omega \right) +\beta
^{2}\int\limits_{\Omega }\left\vert z_{1}\right\vert ^{2}d\Omega \\
&& \\
&\leq &\frac{2\epsilon ^{\ast }}{1-2\epsilon _{2}}\left\Vert \Phi
\right\Vert _{\mathbf{H}}^{2}+C_{\epsilon _{2}^{\ast }}\left[ \frac{1}{%
\left\vert \beta \right\vert ^{\frac{2-4\delta }{1+2\delta }}}\left\Vert
z_{1}\right\Vert _{H^{1}(\Omega )}^{2}+\left\vert \beta \right\vert
^{8}\left\Vert \theta \right\Vert _{H^{1}(\Gamma _{0})}^{2}+\left\vert \beta
\right\vert ^{8}\left\Vert \Phi ^{\ast }\right\Vert _{\mathbf{H}}^{2}\right]
,
\end{eqnarray*}

for $0<\epsilon _{2}<1/2$, and where $C_{\epsilon _{2}}^{\ast }=2C_{\epsilon
_{2}}/(1-2\epsilon _{2})$. If we now take to be $\left\vert \beta
\right\vert $ is sufficiently large; in particular, if 
\begin{equation}
\left\vert \beta \right\vert \geq \mathfrak{B}\equiv \max \left\{ 1,\left(
2C_{\epsilon _{2}}^{\ast }\right) ^{\frac{1+2\delta }{2-4\delta }}\right\} ,
\label{big}
\end{equation}%
we then have 
\begin{equation}
\frac{1}{2}\left( \left\Vert z_{1}\right\Vert _{H^{1}(\Omega
)}^{2}+\left\Vert \beta z_{1}\right\Vert _{\Omega }^{2}\right) \leq \frac{%
2\epsilon ^{\ast }}{1-2\epsilon _{2}}\left\Vert \Phi \right\Vert _{\mathbf{H}%
}^{2}+C_{\epsilon _{2}}\left[ \left\vert \beta \right\vert ^{8}\left\Vert
\theta \right\Vert _{H^{1}(\Gamma _{0})}^{2}+\left\vert \beta \right\vert
^{8}\left\Vert \Phi ^{\ast }\right\Vert _{\mathbf{H}}^{2}\right] .
\label{z_1}
\end{equation}

\bigskip

We can use the wave estimate (\ref{z_1}) to in turn refine the inequality (%
\ref{57.5}) for the mechanical displacement: Combining (\ref{57.5}) and (\ref%
{z_1}) with $\left\vert ab\right\vert \leq \delta a^{2}+C_{\delta }b^{2}$
(and rescaling $\epsilon ^{\ast }>0$), we have for $\left\vert \beta
\right\vert \geq \mathfrak{B}$, as given in (\ref{big}),%
\begin{equation}
\left\Vert \mathring{A}^{\frac{1}{2}}w_{1}\right\Vert _{\Gamma _{0}}^{2}\leq
\epsilon ^{\ast }\left\Vert \Phi \right\Vert _{\mathbf{H}}^{2}+C\left[
\left\vert \beta \right\vert ^{8}\left\Vert \theta \right\Vert
_{H^{1}(\Gamma _{0})}^{2}+\left\vert \beta \right\vert ^{8}\left\Vert \Phi
^{\ast }\right\Vert _{\mathbf{H}}^{2}\right] .  \label{mecah}
\end{equation}
In addition, invoking the resolvent relations (\ref{42}), (\ref{44}), and
estimate (\ref{e5}), and (\ref{z_1}) (and rescaling $\epsilon ^{\ast }>0$),
we have for $\left\vert \beta \right\vert \geq \mathfrak{B}$, as given in (%
\ref{big}),%
\begin{eqnarray}
\left\Vert z_{2}\right\Vert _{\Omega }^{2} &\leq &\epsilon ^{\ast
}\left\Vert \Phi \right\Vert _{\mathbf{H}}^{2}+C\left[ \left\vert \beta
\right\vert ^{8}\left\Vert \theta \right\Vert _{H^{1}(\Gamma
_{0})}^{2}+\left\vert \beta \right\vert ^{8}\left\Vert \Phi ^{\ast
}\right\Vert _{\mathbf{H}}^{2}\right] .  \label{z_2} \\
\left\Vert P_{\gamma }^{\frac{1}{2}}w_{2}\right\Vert _{\Gamma _{0}}^{2}
&\leq &C\left( \beta ^{2}\left\Vert \theta \right\Vert _{H^{1}(\Gamma
_{0})}^{2}+\left\Vert \Phi ^{\ast }\right\Vert _{\mathbf{H}}^{2}\right) .
\label{w_2}
\end{eqnarray}
Finally: combining (\ref{teta}) (and applying Cauchy-Schwartz thereto), (\ref%
{mecah}), (\ref{z_1}), (\ref{z_2}), and (\ref{w_2}), and rescaling $\epsilon
^{\ast }>0$, we have for $\left\vert \beta \right\vert \geq \mathfrak{B}$,
as given in (\ref{big}),%
\begin{equation*}
\left\Vert \Phi \right\Vert _{\mathbf{H}}^{2}\leq \epsilon ^{\ast
}\left\Vert \Phi \right\Vert _{\mathbf{H}}^{2}+C\left[ \left\vert \beta
\right\vert ^{8}\left\Vert \theta \right\Vert _{H^{1}(\Gamma
_{0})}^{2}+\left\vert \beta \right\vert ^{8}\left\Vert \Phi ^{\ast
}\right\Vert _{\mathbf{H}}^{2}\right] .
\end{equation*}%
After applying (\ref{teta}) once more, we have for $\left\vert \beta
\right\vert \geq \mathfrak{B}$ as given in (\ref{big}), 
\begin{equation}
\left\Vert \Phi \right\Vert _{\mathbf{H}}^{2}\leq 2\epsilon ^{\ast
}\left\Vert \Phi \right\Vert _{\mathbf{H}}^{2}+C\left\vert \beta \right\vert
^{16}\left\Vert \Phi ^{\ast }\right\Vert _{\mathbf{H}}^{2}.  \label{final}
\end{equation}
Taking $0<\epsilon ^{\ast }<1/4$, we have at last
\begin{equation*}
\left\Vert \Phi \right\Vert _{\mathbf{H}}^{2}\leq C\left\vert \beta \right\vert
^{16}\left\Vert \Phi ^{\ast }\right\Vert _{\mathbf{H}}^{2}.  \label{final}
\end{equation*}
Combining this estimate with Theorem 10 completes the proof of Theorem \ref{US}. \
\ \ $\square $

\end{document}